\def\be{\begin{equation}}
\def\ee{\end{equation}}

\def\bea{\begin{eqnarray}}
\def\eea{\end{eqnarray}}
\def\beastar{\begin{eqnarray*}}
\def\eeastar{\end{eqnarray*}}

\def\V#1{\buildrel #1\over V}
\def\L#1{\buildrel #1\over L}
\def\P#1{\buildrel #1\over P}
\def\S#1{\buildrel #1\over S}
\def\sig#1{\buildrel #1\over \sigma}
\def\Q#1{\buildrel #1\over Q}

\def\Li{{\cal L}}
\def\Ve{{\cal V}}

\def\Se{{\cal S}}

\def\End{\rm End}
\def\Hom{\rm Hom}
\def\Irr{\rm Irr}

\def\ua{\buildrel\alpha \over u}
\def\ub{\buildrel\beta \over u}
\def\uc{\buildrel\gamma \over u}
\def\u0{\buildrel 0\over u}
\def\u#1{\buildrel #1\over u}

\def\ga{\buildrel\alpha \over g}
\def\gb{\buildrel\beta \over g}
\def\gc{\buildrel\gamma \over g}

\def\mua{\buildrel\alpha \over \mu}

\def\mui{\buildrel{\alpha_i} \over \mu}

\def\muu#1{\buildrel #1\over \mu}

\def\qtr{\mbox{tr}_{q}}
\def\qTr{\mbox{Tr}_{q}}

\def\Wa{\buildrel\alpha \over W}
\def\Wb{\buildrel\beta \over W}
\def\Wc{\buildrel\gamma \over W}
\def\Wf#1{\buildrel #1 \over W}

\def\yi{\buildrel i\over y}
\def\yj{\buildrel j\over y}
\def\zi{\buildrel i\over z}
\def\zj{\buildrel j\over z}
\def\z#1{\buildrel #1\over z}

\def\supijl{{\buildrel {(ij)}\over <}_{l}}

\def\cR#1{{\cal R}^{(#1)}}

\def\cU#1{{\cal U}^{(#1)}}

\def\cUS#1{\buildrel #1\over {\cal U}}

\def\Rab{{\buildrel \alpha\beta \over R}\!}

\def\Rtff#1#2{\buildrel #1 #2 \over {\mathaccent 20 {R}}}
\def\Rtmff#1#2{\buildrel #1 #2 \over {\mathaccent 20 {R}}^{-1}}

\def\Rpff#1#2{\buildrel #1#2 \over {R'}}
\def\Rmff#1#2{\buildrel #1#2 \over {R^{-1}}}
\def\ea{\buildrel\alpha \over e}

\def\guf#1{{\buildrel{#1} \over u}}
\def\muf#1{{\buildrel{#1} \over \mu}}

\def\ea{\buildrel\alpha \over e}

\def\id{\mbox{id}}

\newtheorem{definition}{Definition}
\newtheorem{proposition}{Proposition}
\newtheorem{theorem}{Theorem}
\newtheorem{lemma}{Lemma}

\def\bd{\begin{definition}}
\def\ed{\end{definition}}
\def\bp{\begin{proposition}}
\def\ep{\end{proposition}}
\def\proof{{\sl Proof:}\,}
\def\cqfd{$\Box$}

\documentstyle[bezier,psfig]{article}
\title{Link invariants and Combinatorial  Quantization of Hamiltonian Chern
Simons theory.}
\author{E.Buffenoir\thanks{e-mail: buffenoi@orphee.polytechnique.fr},
Ph.Roche\thanks{e-mail: roche@orphee.polytechnique.fr}\cr \cr
 Centre de Physique Theorique Ecole Polytechnique\thanks{Laboratoire Propre du
CNRS UPR 14}\cr
91128 Palaiseau Cedex\cr
France}

\date{\today}
\begin{document}
\maketitle

\begin{abstract}
We define and study the properties of observables associated to any link in
$\Sigma\times {\bf R}$ (where $\Sigma$ is a compact surface) using the
combinatorial quantization of hamiltonian Chern-Simons theory. These
observables are traces of holonomies in a  non commutative Yang-Mills theory
where the gauge symmetry is ensured by a quantum group. We show that these
observables are link invariants taking  values in a non commutative algebra,
the so called Moduli Algebra. When $\Sigma=S^2$ these link invariants are pure
numbers and are equal to Reshetikhin-Turaev link invariants.

Preprint CPTH: RR 367-07-95
\end{abstract}

\section{Introduction}
Since the fundamental discovery by V.Jones in 1984 of a new link invariant,
there has been a tremendous interest and activity
in low dimension topology using field theory  techniques. The original
definition of the Jones Polynomial was
purely combinatorial and a geometrical understanding of it was finally given by
E.Witten in 1989 \cite{W}. He showed that the Jones polynomial could be
interpreted as the correlation function of Wilson loops (i.e traces of
holonomies) in
 Chern-Simons theory. His work opened a new area of research in what is now
called three dimensional topological field theory. Although this theory is
purely topological (i.e in a hamiltonian picture the hamiltonian is zero) and
therefore contains no dynamics, the quantization of this theory is not at all a
trivial task, mainly because there is no direct procedure to quantify this
theory.
The original method of E.Witten is a brilliant use of  path integrals,
heuristic regularization (by a framing)  of Wilson loops  and relations with
conformal field theory.
Although very appealing and having far reaching consequences, his formalism is
not at all mathematically well defined and this is one of the reasons why many
researchers in this field have used other approaches.  These methods can be
roughly divided in two classes: perturbative and non-perturbative methods.
On the one hand perturbative methods have continuously attracted interest
\cite{GMM} and have provided many interesting recent results: generalization of
Gauss invariants, connections with Vassiliev invariants \cite{AxS} etc...
On the other hand, the geometrical quantization program and combinatorial
methods are the main approaches  to quantize nonperturbatively   Chern-Simons
theory.
Combinatorial methods, introduced by Reshetikhin-Turaev \cite{RT1} and
Turaev-Viro \cite{TV} give explicit representations of abstract amplitudes
satisfying algebraic   relations of a topological field theory. The essential
ingredient in their approach is the representation theory of modular Hopf
algebras which provides family of numbers satisfying Yang Baxter equation, 6-j
identities etc...
These combinatorial methods, although completely well defined, are losing
completely the relationship with Chern-Simons theory.

In our  work we will continue the study of a different type of quantization,
named combinatorial quantization of Hamiltonian Chern Simons theory, which has
been introduced by V.V.Fock and A.A.Rosly \cite{FR} and further developped by
A.Y.Alekseev and al in \cite{AGS1,AGS2,AS} and ourself \cite{BR}. This
quantization can be thought of as being a lattice regularization of Chern
Simons theory in the spirit of Wilson. After quantization, gauge invariance is
replaced by gauge invariance under a quantum group, and the lattice variables,
group elements before quantization, are replaced by elements of a non
commutative algebra. The final and central object of our study is a two
dimensional non commutative Yang Mills theory. Elements of this  program was
already described in the abelian case in \cite{ES}.

In section 2 of this work we give a summary of works  on non commutative two
dimensional Yang Mills theory. We associate to each compact triangulated
surface $\Sigma$ a lattice gauge theory which is covariant under a quantum
gauge group.
The algebra $\Lambda$ of gauge fields is non commutative because matrix
elements of gauge fields associated to arbitrary edges are non commuting.
Locality is however preserved in the sense that matrix elements of gauge fields
associated  to edges having no boundary points are commuting elements. Wilson
loops associated to non self-intersecting loops on the surface are defined and
it is shown that these Wilson loops are gauge invariant elements. A non
commutative analogue of the  Yang Mills action is built following the lines of
A.A.Migdal. In the weak coupling regime, exponential of this action is an
analogue of the Dirac delta function which selects gauge fields with zero
curvature. This theory is therefore a topological field theory and the algebra
of observables
(the Moduli algebra) $\Lambda_{CS}$ of this theory is expected to be a new
description of the algebra of observables of Hamiltonian  Chern Simons theory,
i.e when the three manifold is equal to $\Sigma \times {\bf R}.$

In section 3 we generalize this  construction to the case where the loop is an
arbitrary framed link $L$ in $\Sigma\times [0, 1].$  We obtain observables
associated to these framed links which behave as desired: they are gauge
invariant and are invariant under ambiant isotopy.  As a result we obtain a new
type of ribbon invariants which are not pure numbers but take their values in
the algebra $\Lambda_{CS}.$ This algebra is non commutative except in the case
$\Sigma=S^2$
where it is one dimensional. In that case  the ribbon invariants take their
value in the field ${\bf C}.$

The last part of our work gives the proof that these invariants in the case
where $\Sigma=S^2$ are the Reshetikhin-Turaev invariants.

\section{Summary of works on noncommutative two dimensional Yang Mills theory}
In this section we will make constant use of results obtained in
\cite{AGS1,BR,AGS2}.
Let $\Sigma$ be a compact connected  oriented   surface
  and let ${\cal T}$ be a  triangulation  of $\Sigma.$ Let  us denote by ${\cal
F}$  the
faces of  ${\cal T},$  by $\Li$ the set of oriented edges
 and $\Ve$  the set of points (vertices)
of
this triangulation. If $l$ is an edge, $-l$ will denote
the opposite edge and we have $\{l, -l\}\subset \Li.$

If $l$ is an oriented edge it will be convenient to write $l=xy$ where $y$ is
the departure point
of $l$ and $x$ the end point of $l.$ We will write $y=d(l)$ and $x=e(l).$

Let $A$ be a quasitriangular Hopf algebra such that each finite dimensional
A-module is semisimple. Let $\Irr(A)$  be
 the set of all
equivalency classes of  finite dimensional irreducible representations, in
each of these  classes
$\dot{\alpha}$   we will pick out a particular representative $\alpha.$
Let us denote by $\V{\alpha}$ the
vector space on which acts the representation $\alpha.$

We will denote by $\bar{\alpha}$ (resp.  $\tilde{\alpha}$) the right
 (resp. left) contragredient representation associated to
$\alpha$ acting on $\V{\alpha}\!{}^{\star}$ and defined by:
$\bar{\alpha}={}^t\alpha \circ S $ (resp.
$\tilde{\alpha}={}^t\alpha \circ S^{-1}).$
We will also denote by $0$ the representation of dimension $1$ related to
the counit $\epsilon.$

Let, as usual, denote  $R=\sum_{i} a_i\otimes b_i$  the universal $R$~matrix
 of $A$ and let
us define the invertible element $u$ of $A$ by
 $u=\sum_{i}S(b_i)a_i$ (properties of $u$ can be found in \cite{Dr}.)
Two important elements of $A$ are the ribbon central  element $v$ defined by
$v^2=uS(u)$
 and
the element $\mu=uv^{-1}.$  It will be convenient to define the endomorphism
$\mua=\alpha(\mu)$ and the complex number $v_{\alpha}$ by
$\alpha(v)=v_{\alpha}1.$
If $(\beta)\in Irr(A)^{\times n}$ we will use the notation $\V{(\beta)}$ to
denote the space $\V{(\beta)}=\otimes_{i=1}^{n}\V{\beta_i}$ and
$\muf{(\beta)}=\otimes_{i=1}^n \muf{\beta_i}.$

Let $R'=\sigma (R)$ where $\sigma$ is the permutation operator acting
on $A\otimes A,$ we will use the standard notation:
 \be
R^{(+)}=R, \,R^{(-)}=R'{}^{-1}
\ee

 The $q-$dimension of $\alpha$ is defined by $[d_{\alpha}]=
 tr(\alpha(\mu)).$

  Let $(\ea{}_i\vert i=1...dim \V{\alpha})$ be a particular basis of
 $\V{\alpha},$  and
 $(\ea{}^i \vert i=1...dim \V{\alpha})$ its dual basis.
We will define the linear forms $\ga{}^j_i=<\ea{}^j\vert \alpha(.)\vert
\ea{}_i>.$

The existence of $R$ implies that they satisfy the exchange relations:

\begin{equation}
\Rab_{12} \ga_1 \gb_2 = \gb_2 \ga_1 \Rab_{12},\label{Rgg}
\end{equation}

 also equivalent to:
\begin{equation}
\Rab^{(-)} \ga_1 \gb_2 = \gb_2 \ga_1 \Rab^{(-)}\label{Rgg'},
\end{equation}
where $\Rab=(\alpha\otimes \beta)(R)$ and $\Rab^{(-)}=(\alpha\otimes
\beta)(R^{(-)}).$

Let $\Gamma$ be the restricted dual of $A:$ it is by definition
the  Hopf algebra generated as a vector space by the elements $\ga{}^i_j.$

The action of the coproduct on these elements is:

\begin{equation}
\Delta (\ga{}^i_j)=
 \sum_k \ga{}^i_k \otimes \ga{}^k_j\label{coproduit}.
\end{equation}

$\V{\alpha}$ can be endowed with a structure of right  comodule over $\Gamma:$
\be
\Delta_{\alpha}(\ea_i)=\sum_{j} \ea_j\otimes \ga{}^j_i.
\ee

Let $\alpha, \beta$ be two fixed elements of $\Irr(A),$
by assumption finite dimensional representations are
completely reducible, therefore we can write:
\be
\alpha\otimes\beta=
\bigoplus_{\gamma\in \Irr(A)}N_{\alpha\beta}^{\gamma}\, \gamma .
\ee

Let us define, for each $\gamma,$
$(\psi^{\alpha,\beta}_{\gamma,m})_{m=1...N_{\alpha\beta}^{\gamma}}$
a basis of $\Hom_A(\V{\alpha}
\otimes \V{\beta},\V{\gamma})$
 and  $(\phi^{\gamma,m}_{\alpha,\beta})_{m=1...N_{\alpha\beta}^{\gamma}}$ a
 basis of
$\Hom_A(\V{\gamma},\V{\alpha} \otimes \V{\beta}):$

\begin{equation}
\V{\alpha} \otimes \V{\beta}
 \buildrel {\psi^{\alpha,\beta}_{\gamma,m}} \over \longrightarrow
\V{\gamma}
 \buildrel {\phi^{\gamma,m'}_{\alpha,\beta}} \over \longrightarrow \V{\alpha}
\otimes \V{\beta}.\end{equation}

We have the relation:
\be
\ga_1\gb_2=\sum_{\gamma,m}\phi^{\gamma,m}_{\alpha,\beta}
\gc \psi^{\alpha,\beta}_{\gamma,m}.
\ee
We can always assume that these interwiners satisfy the following relations:
\begin{eqnarray}
\sum_{m,\gamma}
\phi^{\gamma,m}_{\alpha,\beta}\psi^{\alpha,\beta}_{\gamma,m}&=&
id_{\V{\alpha}\otimes\V{\beta}}\label{cg1}\\
\psi^{\alpha,\beta}_{\gamma',m'}\phi^{\gamma,m}_{\alpha,\beta}&=&
id_{\V{\gamma}}\delta_{m'}^m\delta_{\gamma'}^{\gamma}\label{cg2}\\
\phi^{\gamma,m}_{\beta,\alpha} &=& \lambda_{\alpha\beta\gamma}
P_{12}\Rab_{21}^{-1}\phi^{\gamma,m}_{\alpha,\beta} \label{cg3}\\
\psi_{\gamma,m}^{\beta,\alpha} &=& \lambda_{\alpha\beta\gamma}^{-1}
\psi_{\gamma,m}^{\alpha,\beta} \Rab_{21}  P_{12}\label{cg4}\\
\psi_{\gamma,m}^{\beta,\alpha}&=&\sum_{m'}(M^{-1})_m^{m'}{[d_{\gamma}]^{1\over
2}[d_{\beta}]^{1\over 2}\over [d_{\alpha}]^{1\over
2}}(\psi^{\beta,\bar\beta}_{0}\otimes id_{\V{\gamma}})(id_{\V{\beta}}\otimes
\phi^{\alpha,m}_{\bar\beta,\gamma})\label{cg5}\\
\phi^{\gamma,m}_{\beta,\alpha}&=&\sum_{m}M^m_{m'}{[d_{\gamma}]^{1\over
2}[d_{\beta}]^{1\over 2}\over [d_{\alpha}]^{1\over
2}}(id_{\V{\beta}}\otimes\psi^{\bar\beta,\gamma}_{\alpha,m'})
(\phi^{0}_{\beta,\bar\beta}\otimes id_{\V{\gamma}})\label{cg6}
\end{eqnarray}
where $\lambda_{\alpha \beta
\gamma}=(v_{\alpha}v_{\beta}v_{\gamma}^{-1})^{1/2}$ and $M\in GL
(N_{\alpha\beta}^{\gamma}).$

\bd[Gauge symmetry algebra]
Let us define for $z\in \Ve,$ the Hopf algebra  $\Gamma_z=\Gamma\times \{z\}$
 and $\hat
\Gamma=\bigotimes_{z\in \Ve} \Gamma_{z}.$ This Hopf algebra  is
 ``the gauge
symmetry algebra.''
\ed
If $z$ is a vertex we shall write $\ga_z$ to denote the embedding of the
 element $\ga$ in
$\Gamma_z.$

In order to define the non commutative analogue of  algebra of gauge fields we
have to
endow the triangulation with an additional structure \cite{FR}, an order
between edges incident to
each vertex, the {\sl cilium order}.

\bd[Ciliation]
A ciliation of the triangulation is an assignment of a cilium
$c_z$ to each vertex $z$   which consists in a non zero tangent  vector at z.
The orientation of the  surface
defines a canonical cyclic order of the links admitting $z$ as departure or end
point. Let $l_1, l_2$ be links incident to a common vertex $z,$
the strict partial cilium order $<_{c} $ is defined by:

 $l_1<_{c}l_2$ if $l_1\not=l_2, -l_2$ and  the unoriented edges
 $c_z,l_1,l_2$ appear
in the cyclic order defined
by the orientation.
\ed
If $l_1, l_2$ are incident to a same vertex $z$ we define:

$$\epsilon(l_1,l_2)=\left\{ \begin{array}{ll}
+1 &\,\mbox{if}\,l_1<_{c} l_2\\
 -1& \,\mbox{if}\,l_2<_c l_1
\end{array}
\right. $$

\bd[Gauge fields algebra]
The algebra of  gauge fields \cite{AGS1} $\Lambda$ is the algebra
generated by the elements $\ua\!\!(l)^i_j$ with $l\in \Li,\alpha\in
\Irr(A),i,j=1\cdots dim \V{\alpha}$ and satisfying the following determining
relations:

\smallskip

{\bf Commutation rules}
\begin{eqnarray}
& &\ua (xy)_1  \ub (zy)_2 \Rab_{12} =
 \ub (zy)_2 \ua (xy)_1 \label{EE}\\
 & &\ua (xy)_1 {}\Rab_{12}^{-1} \ub (yz)_2 =\ub (yz)_2 \ua (xy)_1\label{ES}\\
& &\Rab_{12} \ua (yx)_1  \ub (yz)_2 =  \ub (yz)_2 \ua (yx)_1\label{SS}\\
& &\, \forall \,\,(yx), (yz) \in \Li\, x\not= z \,\,\,{\rm and}\,\,\,
xy<_{c}yz\nonumber\\
& &\ua(l)\ua(-l)=1\label{ES=1}\\
& &\forall \,\,l \in \Li, \nonumber\\
& &\ua (xy)_1 \ub (zt)_2 = \ub (zt)_2  \ua (xy)_1 \label{Udisjoint}\\
& &\forall\,\, x, y, z, t \mbox{ pairwise distinct in}\, \Ve\nonumber
\end{eqnarray}

{\bf Decomposition rule}
\bea
\ua(l)_1 \ub(l)_2&=&\sum_{\gamma,m}\phi^{\gamma,m}_{\alpha,\beta}
\uc(l)\psi^{\beta,\alpha}_{\gamma,m}
  P_{12}, \label{UCG}\\
\u0(l)&=&1,\,\,\forall l\in \Li.
\eea
\ed

 Gauge covariance of gauge fields comes from the the property that $\Lambda$ is
a right $\hat \Gamma$ algebra comodule defined by the morphism of algebra
 $\Omega:\Lambda\rightarrow \Lambda\otimes\hat\Gamma$ :
\be
\Omega(\ua(xy))=\ga_x \ua(xy) S(\ga_y).
\ee
The subalgebra of gauge coinvariant elements of $\Lambda$ is denoted
 $\Lambda^{inv}.$
Moreover it can be shown that $\ua(-l)=\mua{}^{-1 t}\u{\bar\alpha}(l).$
If $z$ is a vertex we will define $\Omega_z:\Lambda\rightarrow
\Lambda\otimes\Gamma_{z}$ to be equal to $\Omega_{z}=(id\otimes p_z)\Omega,$
where $p_z:\Gamma\rightarrow \Gamma_z$ is the morphism of algebra defined by
$p_z=\otimes_{x\in \Ve, x\not= z} \epsilon_{x}$.

It was shown
(provided some assumption on the existence of a basis of $\Lambda$ of a special
type)
that there exists a unique non zero linear form $h\in \Lambda^{\star}$
satisfying:

\begin{enumerate}
\item (invariance) $(h\otimes id)\Omega(a)= h(a)\otimes 1 \,\,\forall a\in
\Lambda$
\item (factorisation) $h(ab)=h(a)h(b)\\
\forall a\in\Lambda_X, \forall b\in\Lambda_Y,\forall X, Y\subset L,\,\,
(X\cup {-X}) \cap (Y\cup {-Y})=\emptyset$
  \end{enumerate}
(we have  used the notation $\Lambda_{X}$ for $X\subset {\cal L}$ to denote the
subalgebra of $\Lambda$ generated as an algebra by $\ua(l)$ with $l\in X$ and
$\alpha \in \Irr(A)$).

It can be evaluated on any element using the formula:
\be
h(\ua(l){}^i_j)=\delta_{\alpha,0}
\ee
where $0$ denotes the trivial representation of dimension $1,$ corresponding to
the counit.

It is convenient to use the notation $\int d h$ instead of $h.$
The following formula is quite important :
\be
h(\ua(l)_1\mua_2 \ua(-l)_2)={1\over [d_{\alpha}]}P_{12}.
\label{ortho}
\ee
We will use this linear form $h$  in section 4 to compute link invariants.

A path $P$ (resp. a loop $P$) is a connected path (resp. a loop) in the graph
attached
 to the triangulation of $\Sigma$, it will also denote equivalently the
continuous curve (resp. loop) in $\Sigma$ defined by the links of $P.$
Following the definition for links, the departure point of $P$ is denoted
$d(P)$ and its
endpoint $e(P).$
A colored path is a couple $(P, \alpha)$ where $P$ is  a path and $\alpha$ is
an element of $\Irr(A).$ In the rest of this work, we will use as a shortcut
the word path instead of colored path. This should cause no confusion.

Properties of path and loops such as self intersections, transverse
intersections will always be understood as properties satisfied by the
corresponding curves on $\Sigma.$

Let $x_0,\cdots x_n$ be points of $\Ve$ such that $x_{i+1}x_i$ is an edge of
the triangulation, this collection of points defines a path  $P=[x_n,...,x_0],$
 with departure point $x_0$ and end point $x_{n}$. In \cite{BR} we defined
 the sign $\epsilon(x_i, P)=\epsilon((x_{i+1}x_{i}),(x_{i}x_{i-1})).$

If $P$ is a simple  path $P=[x_n,\cdots,x_0]$ with $x_0\not= x_n$, we can
define the holonomy along $P$ by
 \be
\ua_P=v_{\alpha}^{{1\over 2}\sum_{i=1}^{n-1}\epsilon(x_i, P)}\prod_{p=n}^1
\ua(x_px_{p-1}).
\ee

When $C$ is a simple loop $C=[ x_{n+1}=x_0,x_n,\cdots,x_0],$
we will define  the holonomy along $C$  by
 \be
\ua_C=v_{\alpha}^{{1\over 2}(\sum_{i=1}^n \epsilon(x_i, C)-
\epsilon(x_0,C))}\prod_{p=n}^{0} \ua(x_{p+1}x_{p}).
\ee

In \cite{BR} we defined an element of $\Lambda$, which we called {\sl Wilson
loop} attached to $C:$

\be
\Wa_C= tr(\mua \ua_C).\label{Wilsonloop1}
\ee

This element is gauge invariant and moreover
it does not depend on the departure point of  the loop $C.$
This last property can be easily shown using another equivalent expression of
$\Wa_C:$

\be
\Wa_{C}=\omega_{\alpha}(C)
tr_{\V{\alpha}^{\otimes {n+1}}}(\mua{}^{\otimes {n+1}}\prod_{i=n}^1 P_{ii-1}
(\prod_{j=n}^{1}\ua(x_{j+1} x_{j})_j A_j)
 \ua(x_1 x_0)_{0})\label{Wilsonloop2}
\ee
where $A_{j}$ is the matrix $R_{jj-1}^{(\epsilon(x_{j},C))-1}$ and
$\omega_{\alpha}(C)=v_{\alpha}^{-{1\over 2}(\sum_{x\in C}\epsilon(x,C))}.$

The equivalence between relations (\ref{Wilsonloop1}) and (\ref{Wilsonloop2})
uses the simple identity:
\be
tr_{1}(\mua_1
P_{12}(\alpha\otimes\alpha)(R^{(\epsilon)-1}))=v_{\alpha}^{\epsilon}
\id_{\V{\alpha}}
\ee
where $\epsilon=\pm1.$

Remark 1: Compared to our first definition of $\omega_{\alpha}(C)$  in
\cite{BR}, we have used a different normalisation, they are related by a simple
 factor $v_{\alpha}.$ The normalisation of the Wilson loops is discussed in
\cite{BR} in great details.

Remark 2: Expression of the type (\ref{Wilsonloop2})  is reminiscent of the
formulas of \cite{Ma} for the conserved charges in the context of quantum lax
pairs.

It can be  shown  that $\Wa_C$ satisfies the following fusion relation:
\be
\Wa_C \Wb_C=\sum_{\gamma\in \Irr(A)} N_{\alpha\beta}^{\gamma} \Wc_C
\label{fusion}
\ee
where $C$ is any simple loop.

It was also shown that the following commutation relations hold
\be
[\Wa_C, \Wb_{C'}]=0 \label{WW}
\ee
when $C, C'$ are simple loops without transverse intersections.

The geometrical content of this last result is very natural and explained in
the sequel.

\smallskip

Remark: we can define the algebra $\Lambda$ for any type of graph provided that
the graph is endowed with a total order  of the link incident to each vertex.
In particular we can consider any triangulation of any manifold of dimension
greater than three, and define a non commutative lattice gauge theory
associated to it. Unfortunately in that case  we lose the property (\ref{WW})
which is of central importance to define a non commutative Yang Mills action
commuting with gauge invariant elements. This is the reason which prevent us to
extend the present formalism to higher dimensions.

\smallskip

Although the structure of the algebra $\Lambda$ depends on the ciliation, it
has been shown in \cite{AGS1} that the algebra $\Lambda^{inv}$ does not depend
on it up to isomorphism. This is completely consistent with the approach of
V.V.Fock and A.A.Rosly: in their work the graph needs to be endowed with a
structure of ciliated fat graph in order to put on the space of graph
connections ${\cal A}^l$ a structure of Poisson algebra compatible with the
action of the
gauge group $G^{l}.$ However, as a Poisson algebra ${\cal A}^l/G^{l}$ is
canonically isomorphic to the space ${\cal M}^G$ of flat connections modulo the
gauge group, the  Poisson structure of the latter being independent of any
choice of r-matrix \cite{FR}.

In \cite{BR} we introduced a Boltzmann weight attached to any simple loop $C$
and defined by:
\be
\delta_{C}=\sum_{\alpha\in Irr(A)}[d_{\alpha}] \Wa_{C}.
\ee
It was shown that this element satisfies a delta function  property:
\begin{eqnarray}
\delta_{C}\ua_{C}\!{}^i_j&=&\delta_{ij}\delta_{C}.\label{delta}
\end{eqnarray}

We were led to define an element $a_{YM},$ generalizing to our setting the
exponential of the Yang Mills action in the topological limit and  defined by:
$$a_{YM}=\prod_{f\in {\cal F}}\delta_{\partial f}$$
(note that from the relation (\ref{WW}) the elements of this product are
pairwise commuting).

This element satisfies the equation:
\be
a_{YM}\ua_{C}=1 a_{YM},\label{flat}
\ee
for each homologically trivial simple loop on $\Sigma.$
This element is the non commutative analogue of the projector on the space of
flat connections.
The argument leading to commutation relation (\ref{WW}) can be generalized, and
it was proved in \cite{AGS2}
that $\delta_{\partial f}$ for $f\in {\cal F}$ is a central element of
$\Lambda^{inv}.$

The algebra $\Lambda_{CS}=\Lambda^{inv}a_{YM}$ was shown \cite{AGS2} to be
independent , up to isomorphism, of the triangulation. As  a result it was
advocated that $\Lambda_{CS}$ is the algebra of observables of the
Chern Simons theory on the manifold $\Sigma\times {\bf R}.$ This is supported
by the topological invariance of $\Lambda_{CS}$ (i.e this algebra depends only
on the topological structure of the surface $\Sigma$ and not on the
triangulation) and the flatness of the connection.

This geometrical representation of $\Lambda_{CS}$ is particularly appealing.
In particular the element $\Wa_{C}$ that we already  built should be
interpreted as
 being the observable associated to  Wilson loop of horizontal curves,
i.e  loops in $\Sigma\times \{t \}.$
The time $t$, which is the  third coordinate in $\Sigma\times {\bf R},$
manifests itself in the algebraic  point of view as an element used to order
the observables:
if $C_{1},\cdots, C_n$ are colored loops on $\Sigma,$  the element
$W_{C_1}\cdots W_{C_n}$ is the observable associated to the link
 $L=\cup_{i=1}^n\{(C_i,t_i)\}$ where $t_1<\cdots <t_n.$
Note that we are free to choose any time $t_i$  provided that they respect the
order $t_1<\cdots <t_n,$ this relative independance on the time variable is a
simple consequence of the vanishing of the hamiltonian of Chern-Simons theory.

In particular, if  $C$ and $C'$ are curves on $\Sigma$ with no intersection
points, the curve $(C,t)$ and $(C',t')$ never intersect, as a result we obtain
that

\be
\Wa_C \Wb_{C'}= \Wb_{C'}\Wa_C
\ee
 which is the result (\ref{WW}) (note that this last result  was proved in the
more general case where there is no transverse intersections).

Our aim is now to construct in the algebra $\Lambda_{CS}$ the observables
related to  Wilson loop associated to any link in $\Sigma\times [0,1].$

\smallskip

Remarks: 1. In order to simplify our work we have assumed that the surface has
no punctures. This situation can be handled using vertical lines as shown in
\cite{AGS2}.

2. The definition of the element $\delta_{C}$ assumes that the algebra $A$ has
a finite number of irreducible representations. Unfortunately we want to apply
our formalism to the case where $A={\cal U}_q({\cal G}).$  When $q$ is generic
we can however formally bypass this technical problem using the formal
properties of $\delta_{C}$ such as (\ref{delta}). The only infinity which can
possibly occur comes from the square of $\delta_C$ due to the relation:
\be
\delta_{C}^2=(\sum_{\alpha\in Irr(A)} [d_{\alpha}]^{2})\delta_{C}
\ee.

The only sensible way to cure this problem seems to work with $q$ being a root
of unity and to truncate the spectrum either by quotienting by an appropriate
ideal or by using the formalism of weak quasi Hopf algebra as shown in
\cite{AGS1,AGS2}.

\section{Construction of observables $W_{L}$ associated to links $L$ in
$\Sigma\times [0,1]$}
\subsection{Links and chord diagrams}

A link in $\Sigma\times [0,1]$ is an embedding of $(S^1)^{\cup p}$ into
$\Sigma\times [0,1].$
On the set of links we can define a composition law \cite{T}, denoted $*$
defined as follows: let $j_{[a,b]}$ be any increasing diffeomorphism from
$[0,1]$ to $[a,b],$ and let us denote  ${\bar j}_{[a,b]}=id_{\Sigma}\times
j_{[a,b]},$ the composition $L*L'$ is defined by
\be
L*L'={\bar j}_{[0,{1\over 2}]}(L)\cup{\bar j}_{[{1\over 2},{1}]}(L')
\ee

which is an element of
$$\Sigma\times [0,{1\over 2}]\cup_{\Sigma\times\{{1\over 2}\}}\Sigma\times
[{1\over 2},1]=\Sigma\times [0,1].
$$
When the links are considered up to ambiant isotopy this composition is
associative and admit the empty link as unit element. This composition law is
commutative if and only if $\Sigma$ is homeomorphic to the sphere.

In the sequel we will use as a shortcut the term link to  denote  a colored
link in $\Sigma\times [0,1]$ (a link in $\Sigma\times [0,1]$ with an element of
$\Irr(A)$ attached to each connected components of $L$) such that the
projection
 of $L$ on $\Sigma$ is a union of loops on $\Sigma$ in generic position (i.e no
more than double points and
 transverse intersections at these points). Let us denote by
$(\L{i})_{i=1\cdots p}$ the
 connected components of the link $L,$  $\alpha_i\in \Irr(A)$ the color of this
component, and by $\P{i}$ the colored loop
obtained by projecting $\L{i}$ on $\Sigma.$ It is very convenient to associate
to the link $L$ a colored chord diagram \cite{Va,BN,AMR1} which will encode
intersections of
the loops $\P{i}.$ This chord diagram is constructed as follows: the projection
of the link on $\Sigma$ defines $p$ colored loops on $\Sigma$ with transverse
intersections, this configuration of loops defines uniquely a colored chord
diagram by the standard construction. Let us denote by $(\S{i})_{i=1\cdots p}$
the coloured circles of the chord diagram. Each circle $\S{i}$ is oriented, we
will denote by $(\yi_j)_{j=1\cdots n_i}$ the intersection points  of the
circle $\S{i}$ with the chords. We will assume that they are labelled in such a
way that $\yi_j$ appears before $\yi_{j+1}$ with respect to the cyclic order
defined by the orientation of the circles.
Let $Y=\cup_{i=1}^p \{\yi_j,j=1\cdots n_i\}$, we define a relation $\sim$ on
the set $Y$
by :
\be
y\sim y' \,\,\mbox{if and only if} \,y\,\mbox{and}\, y' \,\mbox{are  connected
by a chord.}
\ee
We will denote by $\varphi$ the
immersion of the chord diagram in $\Sigma$, in particular we have
$\P{i}=\varphi(\S{i}).$ Each intersection point of the projection of $L$ on
$\Sigma$ has exactly two inverse images by $\varphi$ in the chord diagram and
these points are  linked by a unique chord. If $p, q$ are two points of $\S{i}$
we will use the notation $[pq]$ to denote the oriented arc segment of  $\S{i}$
having $q$ as departure point and $p$ as endpoint. In the rest of this work we
will assume that $\varphi[\yi_{j+1}\yi_j]$ contains at least two links, for
all $i,j.$ This allows us to find a point $\zi_j\in \S{i}$ such that
$\zi_j\in ]\yi_{j}\yi_{j-1}[$ and $\varphi(\zi_j)$ is a vertex of the
triangulation. This is a purely technical assumption which could be easily
removed.
Each circle $\S{i}$ is the union of $2n_i$ oriented arc segments of type
$[\yi_{j}\zi_j]$ and $[\zi_j\yi_{j-1}],$ let us denote by $\Se_i$  this family
of segments,  $\Se=\cup_{i=1}^p\Se_i$ and  $Z=\cup_{i=1}^p \{\zi_j,j=1\cdots
n_i\}.$

\medskip

\centerline{\psfig{figure=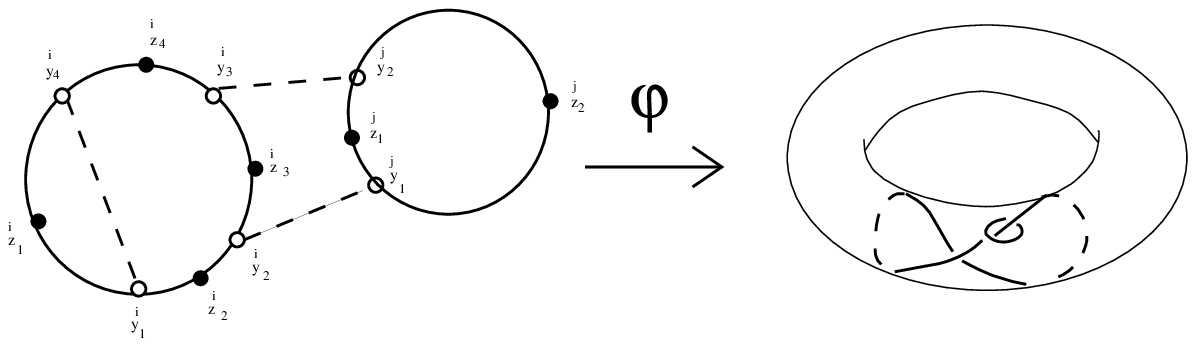}}

\smallskip

\centerline{Fig 1}

\medskip

To each segment $s=[pq]$ of $\Se_i$ we will as usual denote its end point
$e(s)=p$ , its departure point $d(s)=q$ and associate to $s$ the vector spaces
$V_{q^{-}}$
and $V_{p^{+}}$ such that $V_{q^{-}}=
V_{p^{+}}=\V{\alpha_i}.$

If $a$ is a point of  $Y\cup Z$, we define $s(a^{+})$ (resp $s(a^{-})$) to be
the unique element of $\Se$ such that $e(s(a^{+}))=a$ (resp $d(s(a^{-}))=a).$
We shall also use the arc segment $s(a)=s(a^{+})\cup s(a^{-}).$

Let $a\in Y\cup Z,$  and $\xi\in \{+,-\},$ we define $l(a^{\xi})$ to be the
link of the path $\phi(s(a^{\xi}))$ incident to $\phi(a).$

$\Se$ being a finite set, we  choose on it a total ordering and denote by
$\triangleleft$ the strict ordering associated to it.
This ordering allows us to define two vector spaces $V_-$ and $V_+:$
\be
V_-=\bigotimes_{s\in\Se}V_{d(s)^-}=\bigotimes_{x\in Y\cup
Z}V_{x^-}\,\mbox{and}\,
V_+=\bigotimes_{s\in\Se}V_{d(s)^+}=\bigotimes_{x\in Y\cup Z}V_{x^+},
\ee
where the order in the tensor product is taken relativ to $\triangleleft.$

Let $a, b\in Y\cup Z$ and $\xi, \eta\in \{+, -\},$ and assume that
$\phi(a)=\phi(b),$ we will use as a shortcut the notation:
\be
\epsilon(a^{\xi}b^{\eta})=\epsilon(l(a^{\xi}), l(b^{\eta})).
\ee

We define the space  $\Lambda_{\Se}$ by :
\be
\Lambda_{\Se}=\Lambda\otimes \bigotimes_{s\in\Se}
\Hom(V_{d(s)^{-}},V_{e(s)^{+}})
\ee
where the order in the tensor product is taken with respect to $\triangleleft.$

Let $a, b\in Y\cup Z$ and $\xi, \eta\in \{+, -\},$
we shall always use the notation $P_{a^{\xi}b^{\eta}}$ to denote the
permutation operator exchanging the vector spaces $V_{a^{\xi}}$ and
$V_{b^{\eta}}$ in a tensor product of vector spaces.

If $s$ is an element of $\Se_i$ we denote by $j_s$ the canonical injection
$j_s:\Lambda\otimes\Hom (V_{d(s)^{-}},V_{e(s)^{+}})\hookrightarrow
\Lambda_{\Se}.$

Let us define two types of holonomy along $s$:
\begin{itemize}

\item $u_{s}\in\Lambda\otimes\Hom (V_{d(s)^{-}},V_{e(s)^{+}})$ is defined by
$u_{ s}=u_{\varphi(s)},$ (the right handside has already been defined so that
there is no risk of confusion).

\item $U_{s}\in\Lambda_{\Se}$ is defined by $U_{s}=j_{s}(u_{\varphi(s)}).$
\end{itemize}

We have to introduce both of these holonomies because in the rest of this work
some of the constructions are easier to formulate with $U_{s}$ whereas some are
easier to work with $u_s.$
These two points of view were already appearing in our previous work. Indeed if
$C$ is a closed path then $W_C$ which is defined by (\ref{Wilsonloop1}) can
also be written as (\ref{Wilsonloop2}).
The expression (\ref{Wilsonloop1}) only uses the holonomies $u(l)$ ( playing
the role of variables $u_s$) whereas
expression (\ref{Wilsonloop2}) uses only the variables $u(l)_j$ (playing the
role of variables $U_s$.) The definitions and notations we are explaining in
this section have been designed to  include automatically the auxiliary space
previously labelled by a number. As a result the holonomy $U_{s}$ is labelled
by a segment $s$ and the auxiliary space is also labelled by the same segment.
This has the advantage to greatly simplify the notations. The price to pay is
that we have to choose a total order on $\Se$ (the order $\triangleleft$).

To formulate the exchange relations satisfied by the $U_{s}$ we will introduce
one more definition:
\bd
Let $a, b$ two points of $Y\cup Z$ such that $a\in \S{i}, b\in \S{j}$  and
$\phi(a)=\phi(b),$ let us define the endomorphism:

$R^{(a^{\xi}b^{\eta})}\in \End(V_{a^{\xi}}\otimes V_{b^{\eta}})$ (resp  $
\End(V_{b^{\eta}}\otimes V_{a^{\xi}})$) if $s(a^{\xi})\triangleleft
s(b^{\eta})$ (resp if $s(b^{\eta})\triangleleft s(a^{\xi}))$ by:

$$R^{(a^{\xi}b^{\eta})}=\left\{ \begin{array}{ll}
(\alpha_i\otimes\alpha_j)(R^{(\epsilon(a^{\xi}b^{\eta}))}) &\mbox{if
$s(a^{\xi})\triangleleft s(b^{\eta})$}\\
 P_{a^{\xi}b^{\eta}}(\alpha_i\otimes
\alpha_j)(R^{(\epsilon(a^{\xi}b^{\eta}))})P_{a^{\xi}b^{\eta}}& \mbox{if
$s(b^{\eta})\triangleleft s(a^{\xi})$}
\end{array}
\right. $$

\ed
{}From the previous definitions, the endomorphism  $R^{(a^{\xi}b^{\eta})}$ acts
on the space $V_{a^{\xi}}\otimes V_{b^{\eta}}$ or $V_{b^{\eta}}\otimes
V_{a^{\xi}}$ depending of the order of $s(a^{\xi})$ and $s(b^{\xi})$ with
respect to $\triangleleft.$ It will be sometimes useful to use the notation
$R^{(a^{\xi}b^{\eta})}_{c^{\rho}d^{\sigma}}$ to denote the element $R$ or
$R'{}^{-1}$ (depending on the order
of $s(a^{\xi})$ and $s(b^{\eta})$ with respect to $<_c$) represented on the
space $V_{c^{\rho}}\otimes V_{d^{\sigma}}$ or $V_{d^{\sigma}}\otimes
V_{c^{\rho}}$ (depending on the order of $s(c^{\rho})$ and $s(d^{\sigma})$ with
respect to $\triangleleft.$ Note that the relation
$R^{(a^{\xi}b^{\eta})}_{c^{\rho}d^{\sigma}}R^{(b^{\eta}a^{\xi})}_{d^{\sigma}c^{\rho}}=1 \otimes 1$ always holds true.

\bp[General exchange relations]
The elements $U_{s}, s\in\Se$ satisfy the following exchange relations in
$\Lambda_{\Se}$:
\begin{eqnarray*}
R^{(y_1^{+}y_2^{+})}U_{[y_1z_1]}U_{[y_2z_2]}&=&U_{[y_2z_2]}U_{[z_1y_1]}\,,\forall [y_1z_1],
[y_2z_2]\in \Se, y_1\sim y_2,\\
U_{[z_1y_1]}R^{(y_1^{-}y_2^{+})-1}U_{[y_2z_2]}&=&U_{[y_2z_2]}U_{[y_1z_1]},\,\forall [y_1z_1],
[y_2z_2]\in \Se, y_1\sim y_2\\
U_{[z_1y_1]}U_{[z_2y_2]}R^{(y_1^{-}y_2^{-})}&=&U_{[z_2y_2]}U_{[z_1y_1]},\,\forall [z_1y_1],[z_2y_2]\in\Se, y_1\sim y_2\\
U_{s_1}R^{(d(s_1)^{-}e(s_2)^{+})-1}U_{s_2}&=&U_{s_2}U_{s_1},\,\forall
s_1,s_2\in \Se, e(s_2)=d(s_1),\varphi (d(s_2))\not= \varphi(e(s_1))\\
U_{[y_1z_1]}R^{(z_1^{-}z_1^{+})-1}U_{[z_1y_2]}&=&
U_{[z_1y_2]}R^{(y_2^{-}y_1^{+})-1}U_{[y_1z_1]},,\,\forall [y_1z_1],[z_1y_2]\in
\Se,y_1\sim y_2\\
U_{s_1}U_{s_2}&=&U_{s_2}U_{s_1},\,\forall s_1,s_2\in \Se, \phi(s_1)\cap \phi
(s_2)=\emptyset
\end{eqnarray*}
\ep
\proof It follows straightforwardly from the definition of $U_{s}$ and the
exchange relations of the gauge fields $u(l).$
\cqfd

Up to now we have not used the information coming from the topology of the link
i.e over and under crossings of the projected paths $\P{i}.$ This will be
encoded in the following definition:
\bd
Let $<$ be any fixed strict total order on  $Y,$ we  define a family
$\{\cR{y}\}_{y\in Y}$ of elements of
$\bigotimes_{s\in\Se} \Hom (V_{d(s)^-},V_{e(s)^{+}})$ as follows:
let $\{y, y'\}$ be any pair of points of
$Y$ such that $y\sim y',$ we can always assume (otherwise we just exchange $y$
and $y'$) that $y< y',$

\be
\cR{y}=\left\{ \begin{array}{ll}
R^{(y^{-}y^{+})-1} \,\mbox{if $\varphi (s(y))$ is above $\varphi (s(y'))$}\\
R^{(y^{+}{y'}^{+})} R^{(y^{-}y^{+})-1} R^{(y^{-}{y'}^{+})-1}\, \mbox{if
$\varphi (s(y))$ is under $\varphi (s(y'))$}
\end{array}
\right.
\ee
\be
\cR{y'}=\left\{ \begin{array}{ll}
R^{({y'}^{-}{y'}^{+})-1} \,\mbox{if $\varphi (s(y))$ is above $\varphi
(s(y'))$}\\
R^{({y'}^{-}{y'}^{+})-1}R^{({y}^{-}{y'}^{-})}R^{({y'}^{-}{y}^{+})}_{{y'}^{-}{y}^{-}}  \,  \mbox{if $\varphi (s(y))$ is under $\varphi (s(y'))$}
\end{array}
\right.
\ee
This definition defines completely the elements $\cR{y}$ for any $y\in Y.$
Similarly if $z$ is an element of $Z$ we will
define $\cR{z}=R^{(z^{-}z^{+})-1}.$
\ed
The reader can legimately find the definition of $\cR{y}$ obscure, this
definition should be clear after reading the next lemmas.

At this point it is very important to understand that $\cR{y}$ depends on
numerous orderings namely:
\begin{itemize}
\item the total order $<$ on $Y$
\item the total order $\triangleleft$ on $\Se$
\item the partial order coming from the ciliation
\item the under and overcrossings of the projection of the link.
\end{itemize}
We will sometimes write $\cR{y}(<)$ to make explicit  the dependance in $<.$

Let us define the cyclic permutation operator:
$$\sig{i}=P_{\zi_{1}{}^{+}\yi_{1}{}^+}\prod_{j=2}^{n_i}
P_{\zi_{1}{}^+\zi_{j}{}^+}P_{\zi_{1}{}^+\yi_{j}{}^+}.$$

We  now have  defined the framework necessary to associate to $L$ an element of
$\Lambda$ denoted $W_{L}$ which generalizes the construction of Wilson loops.
We denote by $<_l$  the strict lexicographic order induced on $Y$ by the
enumeration of the connected components of $L$ and a choice of departure point
for each of these components, i.e
$\yi_{p}<_l \,\yj_{q}$  if and only if  $i<j$ or ($i=j$ and $p>q.$)

Let us denote the holononomy along the circle $\S{i}$  by
\be
\cUS{i}=\omega(\S{i})U_{[\zi_1\yi_{n_i}]}\cR{\yi_{n_i}}(<_l)U_{[\yi_{n_i}\zi_{n_i}]}\cR{\zi_{n_i}}(<_l)U_{[\zi_{n_i}\yi_{n_{i-1}}]}
\cdots U_{[\yi_{1}\zi_{1}]},
\ee
where
$\omega(\S{i})=v_{\alpha_i}^{-{1\over 2}
\sum_{p=1}^{n_i}\epsilon(\phi(\yi_p),\P{i})+\epsilon(\phi(\zi_p),\P{i})}.$

This holonomy generalizes the holonomy we associated in \cite{BR} to simple
loops: the elements  $\cR{y}$  contain all the information  on the relative
crossings of the projection of the link $L.$

Let $y$ be a point of $Y$ such that $y\in \S{i}$, let $[zy]=s(y^{-})$ and
$[yt]=s(y^{+}),$ we can define the holonomy
 \be
{\cal U}^{(y)}=v_{\alpha_i}^{-{1\over
4}(\epsilon(\varphi(z),\P{i})+2\epsilon(\varphi(y),\P{i})+\epsilon(\varphi(t),\P{i}))}U_{[zy]}\cR{y}U_{[yt]}.
\ee
which is such that
\be
\cUS{i}=\prod_{m=n_{i}}^{2}({\cal U}^{(\yi_{m})}\cR{\zi_m})
{\cal U}^{(\yi_{1})}.
\ee

The following  important lemmas describe the commutation relations of the
elements ${\cal U}^{(y)}$ which explains the definition of $\cR{y}.$

\begin{lemma}
Let $y$ be a point of $Y$ in $\S{i}$ and let us denote by $y'$ the point of $Y$
connected by a chord to $y,$ and assume that $y<y'.$
We have the following commutation relations:
\begin{itemize}
\item
If  $s(y)\cap s(y')=\emptyset$ and
$\varphi (s(y))$ is above $\varphi(s(y'))$ then
\be
\cU{y}\cU{y'}=\sum_{r} \beta^{r}_{{y'}^{+}}R^{({y'}^{+}y^{+})}
\cU{y'}R^{({y'}^{-}{y}^{+})-1}\cU{y} R^{({y'}^- y^-)}
S^2(\alpha^r_{y^{-}})\label{comm1}
\ee
where $R^{(y^-{y'}^+)}=\sum_r \alpha^r_{y^{-}}\otimes \beta^{r}_{{y'}^{+}},$

\item if  $s(y)\cap s(y')=\{z\},$  with $z=d(s(y))=e(s(y'))$ and
$\varphi (s(y))$ is above $\varphi(s(y'))$ then:
\be
\cU{y}R^{(z^-z^+)-1}\cU{y'}=\sum_{r} \beta^{r}_{{y'}^{+}}R^{({y'}^{+}y^{+})}
\cU{y'}R^{({y'}^{-}{y}^{+})-1}\cU{y} R^{({y'}^- y^-)}
S^2(\alpha^r_{y^{-}}).\label{comm2}
\ee
\end{itemize}

Let $a,b$ be two distincts points of $Y$ we denote by $\tau_{a,b}<$ the total
order on $Y$ obtained from $<$ by permuting $a$ and $b.$

The first relation  can also be written:
\begin{itemize}
\item
if $y$ and $y'$ are points of $Y$ linked by a chord such that $s(y)\cap
s(y')=\emptyset$ with $y< y'$ and $\varphi(s(y))$ above $\varphi(s(y')):$
\be
\cU{y}(<)\cU{y'}(<)=\sum_{r} \beta^{r}_{{y'}^{+}}   \cU{y'}(\tau_{y,y'}<)
\cU{y}(\tau_{y,y'}<)R^{({y}^- {y'}^+)-1}_{y^-{y'}^-} S^2(\alpha^r_{y^{-}})
\ee
where $R^{(y^-{y'}^+)}=\sum_r \alpha^r_{y^{-}}\otimes \beta^{r}_{{y'}^{+}}.$

\item
if $y$ and $y'$ are points of $Y$ linked by a chord such that $s(y)\cap
s(y')=\emptyset$  with $y< y'$ and $\varphi(s(y))$ under $\varphi(s(y')):$
\be
\cU{y}(<)\cU{y'}(<)=\sum_{r} \beta^{r}_{{y}^{+}}  \cU{y'}(\tau_{y,y'}<)
\cU{y} (\tau_{y,y'}<)S(\alpha^r_{{y'}^{-}})R^{({y'}^- {y}^+)}_{{y'}^-{y}^-} .
\ee
\end{itemize}
(Similarly if   $s(y)\cap s(y')=\{z\}$  with $z=d(s(y))=e(s(y'))$ the last two
equations hold true if one replaces the left handside by
$\cU{y}R^{(z^-z^+)-1}\cU{y'}$).

\end{lemma}

\proof

This lemma is a direct consequence of the following computation:
\begin{eqnarray*}
&& U_{[zy]}R^{(y^-{y}^+)-1} U_{[yt]} U_{[z'y']}R^{({y'}^-{y'}^+)-1}
U_{[y't']}=\\
&=& U_{[yt]} U_{[zy]} U_{[y't']}U_{[z'y']}\\
&=&\sum_i U_{[yt]} U_{[zy]}\beta^{i}_{{y'}^{+}}
R^{(y^-{y'}^+)-1}S^2(\alpha^i_{y^{-}})U_{[y't']}U_{[z'y']}\\
&=&\sum_i\beta^{i}_{{y'}^{+}} U_{[yt]} U_{[y't']}U_{[zy]}
U_{[z'y']}S^2(\alpha^i_{y^{-}})\\
&=&\sum_i\beta^{i}_{{y'}^{+}} R^{({y'}^+{y}^+)} U_{[y't']}U_{[yt]}
U_{[z'y']}U_{[zy]}R^{(y^-{y'}^-)}S^2(\alpha^i_{y^{-}})\\
&=&\sum_i\beta^{i}_{{y'}^{+}} R^{({y'}^+{y}^+)} U_{[y't']}U_{[z'y']}
R^{({y'}^-y^+)-1}U_{[yt]}U_{[zy]}R^{(y^-{y'}^-)}S^2(\alpha^i_{y^{-}})\\
&=&\sum_{i} \beta^{i}_{{y'}^{+}}
R^{({y'}^{+}y^{+})} U_{[z'y']}R^{({y'}^-{y'}^+)^{-1}} U_{[y't']}
R^{({y'}^{-}{y}^{+})-1}
U_{[zy]}R^{({y}^-{y}^+)^{-1}} U_{[yt]}
R^{({y'}^- y^-)} S^2(\alpha^i_{y^{-}}).
\end{eqnarray*}
Equations (\ref{comm1},\ref{comm2}) are a straightforward consequence of the
previous computations and the definition of $\cR{y}.$
\cqfd

\begin{lemma}
Let $y_1,y_2$ be two points of $Y$ not connected by a chord and let $y_1',
y_2'$ be the points of $Y$ such that $y_1\sim y_1', y_2\sim y_2',$
if we assume    $s(y_1)\cap s(y_2)=\emptyset,$
we have the following relation:
\be
\cU{y_1}(<)\cU{y_2}(<)=\cU{y_2}(<)\cU{y_1}(<).
\ee
Moreover if we assume  $y_1< y_2$ and $\{y_1',y_2'\}\cap \{y\in Y, y_1<y<y_2
\}=\emptyset$
 the last relation can also be written:
\be
\cU{y_1}(<)\cU{y_2}(<)=\cU{y_2}(\tau_{y_1,y_2}<)\cU{y_1}(\tau_{y_1,y_2}<).
\ee

Of course if $s(y_1)\cap s(y_2)=\{z\}$ with $z=d(s(y_1))=e(s(y_2))$  the
previous equations hold true if we replace the left handside by:
 $\cU{y_1}(<)R^{(z^{-}z^{+})-1}\cU{y_2}(<).$
\end{lemma}

\proof
The first part follows straightforwardly from the definition of $\cR{y}$ and
the assumption $s(y_1)\cap s(y_2)=\emptyset.$
The second part comes from the fact that the relative order of $y_i$ and $y_i'$
with respect to $<$ is the same as the one with respect to $\tau_{y_1,y_2}<.$
\cqfd

\subsection{Definition and first properties of the observables $W_L$ associated
to a link $L\subset \Sigma\times [0,1]$}

\bd[Generalized Wilson loops]
To each link in $\Sigma\times[0,1]$ satisfying the assumptions of section (3.1)
we associate an element $W_L$ of the algebra $\Lambda$ by the following
procedure:
let us denote by ${\cal W}_{L}$ the element
\be
{\cal W}_{L}=
\mu_{\Se}\prod_{i=1}^p\sig{i} \prod_{i=1}^{p}
\cUS{i};
\ee
where $\mu_{\Se}=\bigotimes_{s\in \Se}\mu_{e(s)^+}.$

{}From the definition of $\sig{i}$ it follows that:
${\cal W}_{L}\in \Lambda\otimes\Hom(V_{-},V_{+}).$

 The Wilson loop associated to the link $L$ is defined by
$W_{L}=tr_{V_{-}}{\cal W}_{L}$ where  $tr_{V_{-}}$ means the partial trace over
the space $V_{-}$ after the natural identification $V_{+}=V_{-}.$
It will be convenient to write
\be
W_L=Tr_q(\prod_{i=1}^{p}\cUS{i})
\ee
\ed

When $L$ is a simple loop in $\Sigma\times \{t\},$ the element $W_{L}$ we just
defined is  equal to the Wilson loop we already defined by equation
(\ref{Wilsonloop1}).

This element satisfies very important properties which are contained in the
following theorem:

\begin{theorem}
Let $L$ be a link in $\Sigma\times [0,1]$ then $W_L$ does not depend
on the labelling of the components nor  on the choice of departure points of
the components. As a result W is a function on the space of links with values
in $\Lambda.$
Moreover this mapping takes its value in
 $\Lambda^{inv}.$

If $L$ and $L'$ are two links, we have the morphism property:
\be
W_{L\star L'}=W_{L}W_{L'}.
\ee
\end{theorem}

\proof

We will first show that   $W_L$ is invariant under relabelling of the connected
components of the link. It suffices to show that $W_L$ is invariant under the
exchange of $\S{i}$ and $\S{i+1}$ for all $i$. Let $j=i+1,$
we can write
\begin{eqnarray*}
W_L&=&\qTr(A\cUS{i}\cUS{i+1}B)\\
&=&\qTr(A\cU{\yi_{n_i}}\cR{\zi_{n_i}}\cdots\cU{\yi_{1}}
\cU{\yj_{n_j}}\cR{\zj_{n_j}}\cdots\cU{\yj_{1}}B).
\end{eqnarray*}
Let $\supijl$ denote the lexicographic order associated to the labelling of the
connected components after the exchange of $\S{i}$ and $\S{i+1},$ it is
obtained from the first lexicographic order $<_l$ by exchanging  the points
$\yi_{n_i}<_l\cdots <_l\yi_{1}$ with the points
$\yj_{n_j}<_l\cdots <_l\yj_1.$ In order to prove invariance under the
permutation of the connected components it suffices to show that:
\begin{eqnarray*}
&& \qTr(A\cU{\yi_{n_i}}(<_l)\cR{\zi_{n_i}}\cdots\cU{\yi_{1}}(<_l)
\cU{\yj_{n_j}}(<_l)\cR{\zj_{n_j}}\cdots\cU{\yj_{1}}(<_l)B)=\\
&&=\qTr(A\cU{\yj_{n_j}}(\supijl)\cR{\zj_{n_j}}\cdots\cU{\yj_{1}}(\supijl)\cU{\yi_{n_i}}
(\supijl)\cR{\zi_{n_i}}\cdots\cU{\yi_{1}}(\supijl)
B).
\end{eqnarray*}
The proof goes as follows:
if $\yi_1$ is not connected to any $\yj_{m}$ then we deduce  from the second
lemma that:
\begin{eqnarray*}
&&\cU{\yi_1}(<_l)\cU{\yj_{n_j}}(<_l)\cR{\zi_{n_j}}\cdots \cU{\yj_{1}}(<_l)=\\
&=&\cU{\yj_{n_j}}(\tau_{\yi_1,\yj_{n_j}}<_l)
\cU{\yi_1}(\tau_{\yi_1,\yj_{n_j}}<_l)\cR{\zj_{n_j}}\cdots \cU{\yj_{1}}(<_l)=\\
&=&\cU{\yj_{n_j}}(\tau_{\yi_1,\yj_{n_j}}<_l)
\cR{\zj_{n_j}}\cU{\yi_1}(\tau_{\yi_1,\yj_{n_j}}<_l)\cdots
\cU{\yj_{1}}(\tau_{\yi_1,\yj_{n_j}}<_l)=\\
&=&
\cU{\yj_{n_j}}(\prod_{p=1}^{n_{j}}\tau_{\yi_1,\yj_{p}}<_l)\cR{\zj_{n_j}}\cdots
\cU{\yj_{1}}(\prod_{p=1}^{n_j}\tau_{\yi_1,\yj_{p}}<_l)
\cU{\yi_{1}}(\prod_{p=1}^{n_j}\tau_{\yi_1,\yj_{p}}<_l).
\end{eqnarray*}

Let $\tau_q$ denote the permutation
$\tau_{q}=\prod_{p=1}^{n_j}\tau_{\yi_q,\yj_{p}},$
we just have proved that if $\yi_1$ is not connected to any $\yj_{m}$ we have:

\begin{eqnarray*}
&& \qTr(A\cU{\yi_{n_i}}\cR{\zi_{n_i}}\cdots\cU{\yi_{1}}
\cU{\yj_{n_j}}\cR{\zj_{n_j}}\cdots\cU{\yj_{1}}B)=\\
&&=\qTr(A\cU{\yj_{n_j}}(\tau_1 <_l)\cR{\zj_{n_j}}\cdots\cU{\yj_{1}}(\tau_1
<_l)\cU{\yi_{n_i}}
(\tau_1 <_l)\cR{\zi_{n_i}}\cdots\cU{\yi_{1}}(\tau_1 <_l)
B).
\end{eqnarray*}
This last equation still holds true even if $\yi_1$ is connected by a chord to
$\yj_{m}.$
We will show it if $m=n_j,$  the other cases are treated with the same method.

Let us assume first that $\varphi(s(\yi_1))$ is above $\varphi(s(\yj_m))$ and
apply  lemma 1:
\begin{eqnarray*}
&& \qTr(A\cU{\yi_{n_i}}\cR{\zi_{n_i}}\cdots\cU{\yi_{1}}
\cU{\yj_{n_j}}\cR{\zj_{n_j}}\cdots\cU{\yj_{1}}B)=\\
&& =\sum_{r}\qTr(A\cU{\yi_{n_i}}\cR{\zi_{n_i}}\cdots
\cR{\zi_2}\beta_{\yj_{n_j}^+}^r\cU{\yj_{n_j}}(\tau_{1}<_l)\cU{\yi_{1}}(\tau_{1}<_l)\times\\
&&
R^{(\yi_{1}{}^-\yj_{n_j}{}^+)}_{\yi_{1}{}^-\yj_{n_j}{}^-}
S^2(\alpha_{\yj_{n_j}{}^-}^r)\cR{\zj_{n_j}}\cdots\cU{\yj_{1}}B)=\\
&& =\sum_{r}\qTr(\beta_{\yj_{n_j}^+}^r A\cU{\yi_{n_i}}\cR{\zi_{n_i}}\cdots
\cR{\zi_2}
\cU{\yj_{n_j}}(\tau_{1}<_l)\cU{\yi_{1}}(\tau_{1}<_l)\times \\
&&
\cR{\zj_{n_j}}\cdots\cU{\yj_{1}}B
R^{(\yi_{1}{}^-\yj_{n_j}{}^+)}_{\yi_{1}{}^-\yj_{n_j}{}^-}
S^2(\alpha_{\yj_{n_j}{}^-}^r))=\\
&&=\qTr(A\cU{\yi_{n_i}}\cR{\zi_{n_i}}\cdots \cR{\zi_2}
\cU{\yj_{n_j}}(\tau_{1}<_l)\cU{\yi_{1}}(\tau_{1}<_l)
\cR{\zj_{n_j}}\cdots\cU{\yj_{1}}B)\\
&&=\qTr(A\cU{\yi_{n_i}}\cR{\zi_{n_i}}\cdots \cR{\zi_2}
\cU{\yj_{n_j}}(\tau_{1}<_l)
\cR{\zj_{n_j}}\cdots\cU{\yj_{1}}(\tau_{1}<_l)\cU{\yi_{1}}(\tau_{1}<_l)B)=\\
&&=\qTr(A\cU{\yi_{n_i}}(\tau_1 <_l)\cR{\zi_{n_i}}\cdots\cU{\yj_{1}}(\tau_1
<_l)\cU{\yi_{n_i}}
(\tau_1 <_l)\cR{\zi_{n_i}}\cdots\cU{\yi_{1}}(\tau_1 <_l)
B).
\end{eqnarray*}
If  $\varphi(s(\yi_1))$ is under $\varphi(s(\yj_m))$ we use the same type of
proof but use instead relation (\ref{comm2}).

Up to now we have shown that we can move $\cU{\yi_{1}}$ to the right with an
exchange of $<_l$ into $\tau_1\!<_l.$ The same arguments apply as well to
$\cU{\yi_{2}}$ up to  $\cU{\yi_{n_i}}.$
We finally end up with the following equation:
\begin{eqnarray*}
&& \qTr(A\cU{\yi_{n_i}}\cR{\zi_{n_i}}\cdots\cU{\yi_{1}}
\cU{\yj_{n_j}}\cR{\zj_{n_j}}\cdots\cU{\yj_{1}}B)=\\
&&=\qTr(A\cU{\yj_{n_j}}(\prod_{p=n_i}^1\tau_p<)\cR{\zj_{n_j}}
\cdots\cU{\yj_{1}}(\prod_{p=n_i}^1\tau_p<)\times \\
&&\cU{\yi_{n_i}}
(\prod_{p=n_i}^1\tau_p<)\cR{\zi_{n_i}}\cdots\cU{\yi_{1}}
(\prod_{p=n_i}^1\tau_p<)
B),
\end{eqnarray*}
which establishes the result because $\supijl=\prod_{p=n_i}^1\tau_p<.$
This ends up the proof that $W_L$ is invariant under relabelling of the
connected components of $L.$

Let us prove now that $W_L$ is invariant under the choice of departure points
of the curve $\P{i}$ used to define $W_L.$
{}From the  first part of the theorem that we just have proved, we can write:
 \be
W_L=\qTr(\cU{\yi_{n_i}}\cR{\zi_{n_i}}\cdots\cU{\yi_{1}}A).
\ee
Once again the idea of the proof is to use  lemmas 1 and 2 to move
$\cU{\yi_{n_i}}$ to the right.

Let us define $c_i<_l$ the strict order deduced from $<_l$ by applying a cyclic
permutation to  $\yi_{n_i},\cdots,\yi_1$ such that
$\yi_{n_{i-1}} c_i<_l\cdots c_i<_l \yi_{1} c_i<_l \yi_{n_i}.$

Assume first that $\yi_{n_i}$ is not connected to any of th points $\yi_p,$
we can therefore write, from lemma 2:
\begin{eqnarray*}
&&\cU{\yi_{n_i}}\cR{\zi_{n_i}}\cU{\yi_{n_{i-1}}}\cdots\cU{\yi_{1}}=\\
&=&\cU{\yi_{n_{i-1}}}\cU{\yi_{n_i}}\cR{\zi_{n_{i-1}}}\cdots\cU{\yi_{1}}\\
&=&\cU{\yi_{n_{i-1}}}\cR{\zi_{n_{i-1}}}\cdots\cU{\yi_{n_i}}\cU{\yi_{1}}\\
&=&\cU{\yi_{n_{i-1}}}\cR{\zi_{n_{i-1}}}\cdots\cU{\yi_{1}}\cR{\zi_{1}}\cU{\yi_{n_i}}
\end{eqnarray*}
which is equal to:
\be
\cU{\yi_{n_{i-1}}}( c_i<_l)\cdots\cU{\yi_{1}}( c_i<_l)\cR{\zi_{1}}
\cU{\yi_{n_i}}( c_i<_l).
\ee
This proof does not work, but the result is still true, if exceptionally
$n_i=2.$ In that case one has to
generalize  lemma 2 to the case where $s(y_1)\cap s(y_2)=\{z,z'\},
 z\not= z'.$ We leave the details to the reader.

Assume now that $\yi_{n_i}$ is  connected to one of the points $\yi_p,$ and
let us choose $p=n_{i-1}$ to show the structure of the proof:
\begin{eqnarray*}
&&W_L=\qTr(\cU{\yi_{n_i}}\cR{\zi_{n_i}}\cdots\cU{\yi_{1}}A)=\\
&=&\qTr(\cU{\yi_{n_{i-1}}}(\tau_{\yi_{n_i}\yi_{n_{i-1}}}<_l)
\cU{\yi_{n_i}}(\tau_{\yi_{n_i}\yi_{n_{i-1}}}<_l)\cR{\zi_{n_{i-1}}}\cdots
\cU{\yi_{1}}A)
\end{eqnarray*}
using lemma 1 and the property of the quantum trace.

Finally we obtain:
\begin{eqnarray*}
&&W_L=\qTr(\cU{\yi_{n_{i-1}}}(\tau_{\yi_{n_i}\yi_{n_{i-1}}}\!\!<_l)
\cR{\zi_{n_{i-1}}}\cdots\cU{\yi_{n_i}}(\tau_{\yi_{n_i}\yi_{n_{i-1}}}\!\!<_l)
\cU{\yi_{1}}(\tau_{\yi_{n_i}\yi_{n_{i-1}}}\!\!<_l)A)\\
&&=\qTr(\cU{\yi_{n_{i-1}}}(\tau_{\yi_{n_i}\yi_{n_{i-1}}}<_l)
\cR{\zi_{n_{i-1}}}\cdots\cU{\yi_{1}}(\tau_{\yi_{n_i}\yi_{n_{i-1}}}\!\!<_l)
\cR{\zi_{1}}\cU{\yi_{n_i}}(\tau_{\yi_{n_i}\yi_{n_{i-1}}}\!\!<_l)A)\\
&&=\qTr(\cU{\yi_{n_{i-1}}}(c_i<_l)
\cR{\zi_{n_{i-1}}}\cdots\cU{\yi_{1}}(c_i<_l)
\cR{\zi_{1}}\cU{\yi_{n_i}}(c_i<_l)A).
\end{eqnarray*}

It is easy to show using the same method that this  result holds as well
whatever the value of  $p$  can be.

We have shown that we can always replace the departure point $\zi_1$ by its
neighboor $\zi_{n_i}$. As a result, by a trivial induction, it proves  that
$W_L$ is invariant under the change of departure point of each connected
components.

\smallskip

Let us  now show that $W_L$ is a gauge invariant element. It is easy to show
that $W_L$ is a co-invariant element under the coaction $\Omega_x$ where
$x\in\Ve \setminus \varphi(Y):$  this is the same standard proof of prop.2 in
\cite {BR}. The only non trivial part comes from the points of $\varphi(Y)$,
i.e the intersection points of the loops $\P{i}.$
Let $y\in Y$, according to the previous part of the proof, we can always write
$W_L=\qtr(\cU{y}A\cU{y'}B)$ where $y\sim y',$ $\varphi(y)$ above $\varphi(y')$
and $A$ and $B$ are elements of the algebra $\Lambda$ containing no variable
$\ua(l)$ with $l$ incident to $\varphi(y).$

We have:
\begin{eqnarray*}
& &\Omega_{\varphi(y)}(W_{L})=\\
&=&\qTr (U_{[zy]} S(g_{y^{-}}) R^{(y^{-}y^{+})-1}g_{y^+}U_{[yt]}A
U_{[z'y']} S(g_{{y'}^{-}}) R^{({y'}^{-}{y'}^{+})-1}g_{{y'}^+}U_{[y't']}B)\\
&=&\qTr (U_{[zy]} g_{y^+} R^{(y^{-}y^{+})-1 }S(g_{y^{-}})U_{[yt]}A
U_{[z'y']}g_{{y'}^+} R^{({y'}^{-}{y'}^{+})-1} S(g_{{y'}^{-}}) U_{[y't']}B)\\
&=&\qTr (g_{y^+} U_{[zy]}  R^{(y^{-}y^{+})-1 } S(g_{y^{-}})g_{{y'}^+} U_{[yt]}A
U_{[z'y']} R^{({y'}^{-}{y'}^{+})-1} S(g_{{y'}^{-}}) U_{[y't']}B)\\
&=&\sum_i\qTr (g_{y^+} U_{[zy]}  R^{(y^{-}y^{+})-1 }
S(g_{y^{-}})\beta_{{y'}^{+}}^i
R^{(y^{-}{y'}^{+})-1}S^2(\alpha_{y^{-}}^i)g_{{y'}^+} U_{[yt]}A\\
&& U_{[z'y']} R^{({y'}^{-}{y'}^{+})-1} S(g_{{y'}^{-}}) U_{[y't']}B)\\
&=&\sum_i\qTr (\beta_{{y'}^{+}}^i g_{y^+} U_{[zy]}  R^{(y^{-}y^{+})-1 }
S(g_{y^{-}})
R^{(y^{-}{y'}^{+})-1}g_{{y'}^+} U_{[yt]}A\\
&& U_{[z'y']} R^{({y'}^{-}{y'}^{+})-1} S(g_{{y'}^{-}}) U_{[y't']}B
S^2(\alpha_{y^{-}}^i))\\
&=&\sum_i \qTr (\alpha_{y^{+}}^i\beta_{{y'}^{+}}^ig_{y^+} U_{[zy]}
R^{(y^{-}y^{+})-1 } g_{{y'}^+}
R^{(y^{-}{y'}^{+})-1} S(g_{y^{-}})U_{[yt]}A\\
&& U_{[z'y']} R^{({y'}^{-}{y'}^{+})-1)} S(g_{{y'}^{-}}) U_{[y't']}B )\\
&=&\qtr (R_{y^+{y'}^+}^{(y^{-}{y'}^{+})-1}g_{y^+} g_{{y'}^+}
U_{[zy]}  R^{(y^{-}y^{+})-1 } R^{(y^{-}{y'}^{+})-1} U_{[yt]}A\\
&& U_{[z'y']} R^{({y'}^{-}{y'}^{+})-1)}  U_{[y't']}B
S(g_{y^{-}})S(g_{{y'}^{-}}))\\
&=&\qTr (g_{{y'}^+}g_{y^+} R^{(y^{-}{y'}^{+})-1}_{y^+{y'}^+}
U_{[zy]}  R^{(y^{-}y^{+})-1 } R^{(y^{-}{y'}^{+})-1} U_{[yt]}A\\
&& U_{[z'y']} R^{({y'}^{-}{y'}^{+})-1}  U_{[y't']}B
S(g_{y^{-}})S(g_{{y'}^{-}}))\\
&=&\qTr ( R^{(y^{-}{y'}^{+})-1}_{y^+{y'}^+}
U_{[zy]}  R^{(y^{-}y^{+})-1 } R^{(y^{-}{y'}^{+})-1} U_{[yt]}A\\
&& U_{[z'y']} R^{({y'}^{-}{y'}^{+})-1}  U_{[y't']}B).
\end{eqnarray*}
As a result we get $\Omega_{\varphi(y)}(W_{L})\in \Lambda\otimes 1.$
Applying $(id\otimes \epsilon)$ and using the comodule definition, it follows
that:
$\Omega_{\varphi(y)}(W_{L})=W_L\otimes 1,$ i.e $W_L$ is an element of
$\Lambda^{inv}.$

Let $L, L'$ two links of $\Sigma\times [0,1],$ and
$\{\S{i}_{L}i=1,\cdots,p_{L}\},$
(resp. $\{\S{i}_{L'} i=1,\cdots,p_{L'}\}$ ) be the circles associated to the
chord diagramm of $L$ (resp. $L'$). The family of circles of the chord diagram
of $L *L'$ is the union $\{\S{i}_{L}i=1,\cdots,p_{L}\}\cup
\{\S{i}_{L'}i=1,\cdots,p_{L'}\}.$ We can choose a labelling of the circles of
the chord diagram of $L* L'$ such that those of $L$ appear before those of
$L'$. Because $\varphi(L)$ is above $\varphi(L')$ the expression of $\cR{y}$
such
 that $\varphi(y)\in \varphi(L)\cap\varphi(L')$  do not connect $L$ and $L'$
i.e
 $\cR{y}=R^{(y^-y^+)-1}.$ From the definition of $W_{L*L'}$  we immediately get
$W_{LL'}=W_L W_{L'}.$
This  fact is the final step in the  proof of the theorem.
\cqfd

\medskip

As already explained in section 2, when $L$ is a simple loop in $\Sigma\times
[0,1],$ $W_L$ has two equivalent expressions: the first one,  called
" expanded form" (eq.\ref{Wilsonloop2}), can be written
as a partial quantum trace over the space
$\Lambda\otimes\bigotimes_{s\in {\cal S}}\Hom(V_{d(s)^-},V_{e(s)^{+}})$ ,
whereas the second one (eq.\ref{Wilsonloop1}),  called "contracted form", is
expressed as a partial quantum trace  over the space
$\Lambda\otimes\bigotimes_{i=1}^p End(\V{\alpha_{i}}).$
Relation between them is a direct consequence of the following result:

if $A\in \Hom(X,Y)$ and $B\in \Hom(Y,Z)$ then $BA=tr_{Y}(P_{Y,Z}B\otimes A).$

The definition of $W_L$ for arbitrary links is naturally defined in terms of
the $U_{s},$ it can of course be written in terms of a quantum partial trace
over $\Lambda\otimes\bigotimes_{i=1}^p End(\V{\alpha_{i}})$, the expression is
quite messy, the only essential point being that  we can switch from one
expression to the other.
This contracted form enjoys the following property:

\begin{lemma}

 Assume that there exists a connected open path  $P\subset \P{i}$   such  that
\begin{enumerate}
\item $e(P)$ and $d(P)$ are not intersection points of the projected link
\item $\forall y \in \varphi^{-1}(P)$ and $y' \in Y\setminus \varphi^{-1}(P)$
with $y \sim y',$  $\phi(s(y))$ is above $\phi(s(y')),$
\end{enumerate}

then the Wilson loop associated to $L$ can be uniquely written as
\be
W_L= tr_{\V{\alpha_i}}(\muf{\alpha_i} v_{\alpha_i}^{-{1\over2}
\epsilon(e(P),\P{i})} u_{P} v_{\alpha_i}^{{1\over2}
\epsilon(d(P),\P{i})} A),
\ee
 where  $A$ is an element of $\Lambda\otimes End( \V{\alpha_i},\V{\alpha_i})$
which components are linear combinations of the matrix element of
$(\ua(l))_{l \in \phi(L)\setminus P}$ with  coefficients
depending only on the set of links in  $\phi(L)\setminus P$ and  the ciliation
at each vertex of  $\phi(L)\setminus P.$
A path satisfying   properties 1) and 2) will be said to be "on  top" of the
link $L.$
\end{lemma}
\proof
{}From the theorem (1), we can always assume that $P\subset \P{1}$, and from
the independance under the choice of the $\phi(z)$ for $z\in Z,$ we can always
assume that $e(P)=\phi(\z{1}_{1}).$
As a result we can write:
$W_{L}=Tr(\mu_{\Se}\prod_{i=1}^p\sig{i}\cU{1}\prod_{j=2}^p \cU{j}).$
Because $P$ is on the top of $L$ we have $\cR{y}=R^{(y^{-}y^{+})-1}$
if $y\in Y$ and $\varphi(y)\in P.$
Using the formula
\be
tr_{V_{e(s_2^{+})}}(\mu_{e(s_2^{+})}
u_{s_1}R^{(d(s_1)^{-}e(s_2)^{+})-1}u_{s_2})=
v_{\alpha_1}^{\epsilon(\phi(e(s_2)),\P{1})}u_{s_1}u_{s_2},
\ee
we obtain the desired result
\cqfd.

\medskip

\bp[Regular isotopy]

The element ${\widehat W}_L=W_L\prod_{F\in {\cal F}}\delta_{\partial F}$ of
$\Lambda_{CS}$ depends only on the regular isotopy class of the link $L,$ i.e
it satisfies the Reidemeister moves of type 0,2,3.
\ep
\proof
Let us first show an already interesting result in itself:
if $L, L'$ are two links In $|sigma\times {\bf R},$ $L_1$ (resp. $L'1$ ) are
connected curves included in $ L$  (resp.$L'$) and  $P$ (resp. $P'$) are the
projections of $L_1$ (resp $L_1'$) such that:
\begin{itemize}
\item $L\setminus L_1=L'\setminus L_1'$
\item $\phi(L)\setminus P=\phi(L')\setminus P'$
\item $e(P)=e(P'),\, d(P)=d(P')$ and none of these points are intersection
 point
\item $P$ (resp. $P'$) is on  top of $L$ (resp. $L'$)
\item the curve $C=PP'^{-1}$ is a closed simple homologically trivial curve,
\end{itemize}
then the following equality is true:
\be
{\widehat W}_L
={\widehat W}_{L'}.
\ee

In order to show this easy result let us denote $e=e(P)=e(P')$ and
$d=d(P)=d(P'),$ and assume that $P$ and $P'$ are colored by $\alpha_i.$

{}From the definition of a path being on top of $L,$ we have:

\begin{eqnarray*}
\prod_{F\in {\cal F}}\delta_{\partial F}W_L&=&
\prod_{F\in {\cal F}}\delta_{\partial F}
tr_{\V{\alpha_i}}( \mui v_{\alpha_i}^{-{1\over2}
\epsilon(e, L)} u_{P} v_{\alpha_i}^{{1\over2}\epsilon(d,{ L})} A)\\
&=&\prod_{F\in {\cal F}}\delta_{\partial F}
v_{\alpha_i}^{-{1\over2}\epsilon(e,{ L})+{1\over2}\epsilon(d,{L})
+{1\over2}\epsilon(e,{ C})-{1\over2}\epsilon(d,{ C})}\times\\
&\times&tr_{ \V{\alpha_i}}( \mui (v_{\alpha_i}^{-{1\over2}
\epsilon(e,{ C})} u_{P}
v_{\alpha_i}^{{1\over2}\epsilon(d,{ C})} u_{P'^{-1}})
u_{P'} A)\\
&=& \prod_{F\in {\cal F}}\delta_{\partial F}
v_{\alpha_i}^{-{1\over2}\epsilon(e,{L})+{1\over2}\epsilon(d,{ L})
+{1\over2}\epsilon(e,{C})-{1\over2}\epsilon(d,{C})}\times\\
&\times& tr_{\otimes_j \V{\alpha_i}}( \mui u_C u_{P'} A)\\
&=& \prod_{F\in {\cal F}}\delta_{\partial F}
v_{\alpha_i}^{-{1\over2}\epsilon(e,{L})+{1\over2}\epsilon(d,{L})
+{1\over2}\epsilon(e,{C})-{1\over2}\epsilon(d,{C})+
{1\over2}\epsilon(e,{{L'}})-{1\over2}\epsilon(d,{{L'}})}\times\\
&\times& tr_{\otimes_j \V{\alpha_i}}(\mui v_{\alpha_i}^{-{1\over2}
\epsilon(e,{{L'}})}u_{P'}
v_{\alpha_i}^{{1\over2}\epsilon(d,{{L'}})} A)\\
&=& \prod_{F\in {\cal F}}\delta_{\partial F} W_{L'}.
\end{eqnarray*}
The last line comes from the identity:
\be
\epsilon(e,C)+\epsilon(e,L')-\epsilon(e,L)=1.
\ee
This ends the proof of the intermediate result.

The proof of Reidemeister moves of type 0,2,3 is a direct byproduct  of the
previous
result, the proof of all these moves being easily  reduced to move
an open strand on the top of the link.
\cqfd.

In \cite{AGS2} it was shown that if ${\cal T}$ and ${\cal T'}$ are two ciliated
triangulations of $\Sigma,$ the algebra $\Lambda_{CS}({\cal T})$ and
$\Lambda_{CS}({\cal T'})$ are isomorphic. Let us denote $i_{{\cal T}, {\cal
T'}}:\Lambda_{CS}({\cal T})\rightarrow \Lambda_{CS}({\cal T'})$  this
isomorphism,
it was shown that this isomorphism preserves the linear form $h$:
if $a\in \Lambda_{CS}(\cal T)$ then $h(a)=h'(i_{{\cal T}, {\cal T'}}(a))$ where
$h$ (resp. $h'$ ) is the linear invariant form on $\Lambda({\cal T})$
(resp. on $\Lambda({\cal T'})$).

\bp[Invariance under the change of triangulation]
Let ${\cal T}$ and ${\cal T'}$ two triangulations of $\Sigma,$ and let
$R({\cal T}, {\cal T'})$  be any triangulation which is a common
refinement of ${\cal T}$ and ${\cal T'}.$
Let $L$ (resp $L'$) be links in $\Sigma$ such that $L$ resp.( $L'$) has all its
projected components composed with links of ${\cal T}$ resp.($ {\cal T'}$).
If $L$ and $L'$ are related by a finite set of moves of type $0,2,3$ involving
the triangulation $R({\cal T}, {\cal T'})$ then:
$i_{{\cal T}, {\cal T'}}({\widehat W}_{L})={\widehat W}_{L'}.$
A direct consequence of this result is the equality of expectation values, i.e:
$h({\widehat W}_{L})=h'({\widehat W}_{L'}).$
\ep
Proof: this is a direct consequence of proposition 2 and the expression of
$i_{{\cal T}, {\cal T'}}.$
\cqfd

 The only move which can not be deduced
from  proposition2 is the first move. Although being reduced to the move of
 strand on top of the link, the strand $P$ which is moved to $P'$ is such that
$e(P), d(P)$ are elements of $\phi(Y)$ or $PP'^{-1}$ is not a simple curve. So
we cannot apply directly the intermediate result.

In the next theorem we study the behaviour of
 ${\widehat W}_L$ under the move of type 1.

\bp[Type I moves]
Let  $L$ be as usual a link in $\Sigma\times [0, 1]$ and $\P{i}$ the set of
projected curves on $\Sigma$ and  let $L^{\propto\pm}$ be another link
 whose projection $\P{i}{}^{\propto\pm}$ differs from $\P{i}$ by a move of type
I (see fig 2) applied to a curve colored by $\alpha_i$,
we have the following relation:
\be
{\widehat W}_{L^{\propto\pm}}=v_{\alpha_i}^{\pm 1} {\widehat W}_{L}.
\ee
\ep

\medskip

\par
\centerline{\psfig{figure=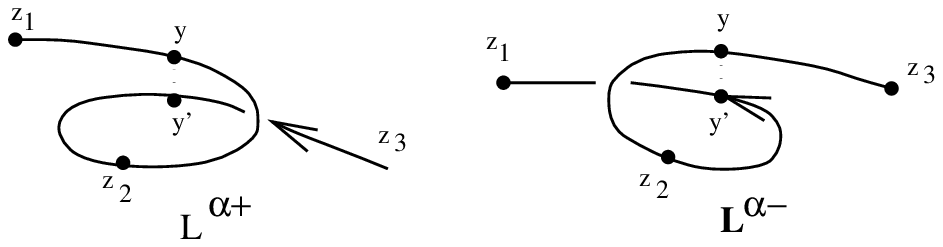}}
\par

\centerline{Fig 2}

\medskip

\proof
Let $C$ denote the closed curve $C=\phi([y' z_2])\cup \phi([z_2 y]),$ we have:

\begin{eqnarray*}
&&\delta_C W_{L^{\propto +}}=\\
&&=\delta_C
v_{\alpha_i}^{-{1\over2}\epsilon({z_1}^{-}{z_1}^{+})+{1\over2}(\epsilon(y^{-}y^{+})
+\epsilon({z_2}^{-}{z_2}^{+})+\epsilon({y'}^{-}{y'}^{+}))+
{1\over2}\epsilon({z_3}^{-}{z_3}^{+})}\times\\
&\times&tr_{ \V{\alpha_i}}( \mui u_{[z_1 y]}u_{[y z_2]}u_{[z_2 y']}u_{[y' z_3]}
A)\\
&=&\delta_C
v_{\alpha_i}^{-{1\over2}\epsilon({z_1}^{-}{z_1}^{+})-{1\over2}(\epsilon(y^{-}y^{+})
+\epsilon({z_2}^{-}{z_2}^{+})+\epsilon({y'}^{-}{y'}^{+}))+
{1\over2}\epsilon({z_3}^{-}{z_3}^{+})}\times\\
&\times&tr_{V_{z_2} \otimes V_{y} \otimes V_{y'} \otimes
\V{\alpha_i}}(P_{z_1^{+} {y'}^{+}}P_{z_1^{+} {z_2}^{+}}P_{z_1^{+} {y}^{+}}
\mua_{y^{+}} \mua_{{z_2}^{+}} \mua_{{y'}^{+}}  \mui \times \\
&\times& u_{[z_1 y]}R^{(y^{-}y^{+})-1}u_{[y
z_2]}R^{({z_2}^{-}{z_2}^{+})-1}u_{[z_2 y']}R^{({y'}^{-}{y'}^{+})-1}u_{[y' z_3]}
A)\\
&=&\delta_C
v_{\alpha_i}^{-{1\over2}\epsilon({z_1}^{-}{z_1}^{+})-{1\over2}(\epsilon(y^{-}y^{+})
+\epsilon({z_2}^{-}{z_2}^{+})+\epsilon({y'}^{-}{y'}^{+}))+
{1\over2}\epsilon({z_3}^{-}{z_3}^{+})}\times\\
&\times& tr_{V_{z_2} \otimes V_{y} \otimes V_{y'} \otimes
\V{\alpha_i}}(P_{z_1^{+} {y'}^{+}}P_{z_1^{+} {z_2}^{+}}P_{z_1^{+} {y}^{+}}
\mua_{y^{+}} \mua_{{z_2}^{+}} \mua_{{y'}^{+}}  \mui \times \\
&\times& u_{[y z_2]}R^{({z_2}^{-}{z_2}^{+})-1}u_{[z_2 y']}u_{[z_1
y]}R^{({y'}^{-}{y}^{-})}R^{({y'}^{-}{y'}^{+})-1}u_{[y' z_3]} A)\\
&=&\delta_C
v_{\alpha_i}^{-{1\over2}\epsilon({z_1}^{-}{z_1}^{+})-{1\over2}(\epsilon(y^{-}y^{+})+\epsilon({y'}^{-}{y'}^{+}))+
{1\over2}\epsilon({z_3}^{-}{z_3}^{+})+{1\over2}\epsilon({y'}^{-}{y}^{+})}\times\\
&\times&tr_{V_{y} \otimes V_{y'}\otimes \V{\alpha_i}}(P_{z_1^{+}
{y'}^{+}}P_{z_1^{+} {y}^{+}} \mua_{y^{+}} \mua_{{y'}^{+}} \mui\times \\
&\times&
v_{\alpha_i}^{-{1\over2}\epsilon({y'}^{-}{y}^{+})}u_{[yz_2]}v_{\alpha_i}^{{1\over2}\epsilon({z_2}^{-}{z_2}^{+})}u_{[z_2 y']}u_{[z_1 y]}R^{({y'}^{-}{y}^{-})}R^{({y'}^{-}{y'}^{+})}u_{[y' z_3]}A)\\
&=&\delta_C
v_{\alpha_i}^{-{1\over2}(\epsilon({z_1}^{-}{z_1}^{+})-{1\over2}(\epsilon(y^{-}y^{+})+\epsilon({y'}^{-}{y'}^{+}))+
{1\over2}\epsilon({z_3}^{-}{z_3}^{+})+{1\over2}\epsilon({y'}^{-}{y}^{+})}\times\\
&\times&tr_{V_{y} \otimes V_{y'}\otimes \V{\alpha_i}}(P_{z_1^{+}
{y'}^{+}}P_{z_1^{+} {y}^{+}} \mua_{y^{+}}  \mua_{{y'}^{+}} \mui \times \\
&\times& \id_{(V_{y^{+}},V_{{y'}^{-}})}
u_{[z_1 y]}R^{({y'}^{-}{y}^{-})}R^{({y'}^{-}{y'}^{+})-1}u_{[y' z_3]} A)=\\
&&(\mbox{after the use of}\,\,
\delta_{C}v_{\alpha_i}^{-{1\over2}\epsilon({y'}^{-}{y}^{+})}u_{[yz_2]}v_{\alpha_i}^{{1\over2}
\epsilon({z_2}^{-}{z_2}^{+})}u_{[z_2
y']}=\delta_{C}\id_{(V_{y^+},V_{y'{}^{-}})}\\
&&\mbox{where $\id_{(V_{y^+},V_{y'^{-}})}$ denote the identity endomorphism
)}\\
&=&\delta_C v_{\alpha_i}^{-{1\over2}(\epsilon({y}^{-}y^{+})
+\epsilon({y'}^{-}{y'}^{+}))-\epsilon({y'}^{-}{y}^{-})+\epsilon({y'}^{-}{y'}^{+})-{1\over2}\epsilon({y}^{-}{y'}^{+})
+{1\over2}\epsilon({y'}^{-}{y}^{+})}\times\\
&\times&tr_{ \V{\alpha_i}}( \mui
v_{\alpha_i}^{-{1\over2}\epsilon({z_1}^{-}{z_1}^{+})}u_{[z_1
y]}v_{\alpha_i}^{{1\over2}\epsilon({y}^{-}{y'}^{+})}u_{[y z_3]}
v_{\alpha_i}^{{1\over2}\epsilon({z_3}^{-}{z_3}^{+})}A)\\
&=& v_{\alpha_i}^{-{1\over2}(\epsilon({y}^{-}y^{+})
+\epsilon({y'}^{-}{y'}^{+}))-\epsilon({y'}^{-}{y}^{-})+
\epsilon({y'}^{-}{y'}^{+})-{1\over2}\epsilon({y}^{-}{y'}^{+})
+{1\over2}\epsilon({y'}^{-}{y}^{+})} \delta_C W_{L}
\end{eqnarray*}
{}From the relation:
$$\epsilon(y'{}^-y'{}^{+})+\epsilon(y^{-}y'{}^-)=\epsilon(y^{-}y'^{+})+1,$$
we can write:
\begin{eqnarray*}
&&-{1\over 2}(\epsilon({y}^{-}y^{+})
+\epsilon({y'}^{-}{y'}^{+}))-\epsilon({y'}^{-}{y}^{-})+
\epsilon({y'}^{-}{y'}^{+})-{1\over2}\epsilon({y}^{-}{y'}^{+})
+{1\over2}\epsilon({y'}^{-}{y}^{+})=\\
&=&{1\over
2}(\epsilon(y'{}^-y^+)+\epsilon(y^+y^-)+\epsilon(y^-y'^{+})+\epsilon(y'{}^+y'{}^-))+1=1,
\end{eqnarray*}
which leads to
$\delta_{C}W_{L}^{\propto +}=v_{\alpha_i} \delta_{C}W_{L}.$
 We finally obtain the relation:
\be
{\widehat W}_{L^{\propto +}}=v_{\alpha_i} {\widehat W}_{L}.
\ee
The proof in the other case is very similar and it can be shown that:
\be
{\widehat W}_{L^{\propto -}}=v_{\alpha_i}^{-1} {\widehat W}_{L}.
\ee
\cqfd

{}From propositions 2 and 3 we immediately obtain that ${\widehat W}_L$ is a
ribbon invariant, the link being endowed with the ``blackboard framing''
associated to the projection $\Sigma\times [0,1]\rightarrow \Sigma.$
Let $n_i$ be the writhe of the ribbon defined by $L_{i},$ Propositions 2 and 3
imply that
\be
I(L)=\prod_{i=1}^p v_{\alpha_{i}}^{-n_{i}} {\widehat W}_{L}
\ee
is a link invariant element of the   algebra $\Lambda_{CS}.$

A particularly simple situation appears when $\Sigma=S^2,$ because in that case
$W_L$ is essentially a number.

Indeed let $\Sigma=S^2,$ then from the independence of the algebra structure of
$\Lambda_{CS}$ under the choice of graph \cite{AGS2} we obtain that
$\Lambda_{CS}(S^2)$ is an algebra of dimension $1,$ (this comes from the fact
that $S^{2}$ is homeomorphic to a disk whose boundary has been identified to a
point). As a result we get
 $\Lambda_{CS}(S^2)={\bf C}\prod_{f\in {\cal T}}\delta_{\partial f}.$
We then deduce that:
\be
{\widehat W}_{L}=w_L\prod_{f\in {\cal T}}\delta_{\partial f},
\ee
with $w_L\in {\bf C}.$

This number $w_L$ satisfies the following equation:
\be
w_{L*L'}=w_{L}w_{L'},
\ee
which is a trivial consequence of the identity $W_{L}W_{L'}=W_{L\star L'}$ and
$dim \Lambda_{CS}(S^2)=1.$

In the next section we will show that $w_L$ is the Reshetikhin-Turaev invariant
(denoted $RT(L)$) of the framed link  $L.$  We will have therefore a new
description of these invariants in term of traces of holonomies of flat
connections in the spirit of the work of E.Witten.

When $\Sigma$ is not homeomorphic to $S^2,$ $\Lambda_{CS}$ is no more a trivial
algebra. It is then a quantization of the space $Fun (Hom (\pi_1(\Sigma), G)/G,
{\bf C}),$ a presentation by generators and relations being given in
\cite{AGS2}. This presentation is simply obtained by generalizing the
definition of $\Lambda$ to arbitrary cell decomposition of the surface. By
taking the simplest one, i.e one 2-cell, 2g edges and one point, these authors
have obtained a nice presentation of the Moduli algebra. This description is
one of the step to obtain the complete set of irreducible unitary
representation of $\Lambda_{CS}$ (even with punctures) as described in
\cite{AS}.
We also have to mention the work \cite{T,AMR2} for other interesting
constructions related to quantization of the moduli space of flat connections.

\section{Relation between the observables $W_L$ with $L\subset S^2\times [0,1]$
and Reshetikhin-Turaev invariants}

The aim of this section is to show that in the case of the sphere the invariant
of link $w_L$ is equal to Reshetikhin-Turaev invariant of the ribbon (endowed
with the blackboard framing) associated to $L$. The idea of the method is very
simple: because ${\widehat W}_L=w_L a_{YM},$ we obtain
$w_L=h(w_L a_{YM})/h(a_{YM}).$
It remains to integrate $w_L a_{YM}$ over all links of the triangulation. This
is quite technical, the final answer being just the Reshetikhin-Turaev
invariant of the link $L$ in  the shadow world.

 Before explaining this proof we will first show that  when  $A={\cal
U}_{q}(sl_{2})$
and $L$ is a link colored with the fundamental representation then
$w_{L}=q^{{3\over 2}\sum_i n_i}P(L, q^2)$  where $P(L,z)$ is the Jones
polynomial of the link $L.$

\bp[Skein relations-Jones Polynomial]
Let  $L_{+}, L_{-}, L_{0}$ be  links in $\Sigma\times [0, 1]$ which coincide
outside a ball and look as in (fig 3) inside the ball.

\medskip

\par
\centerline{\psfig{figure=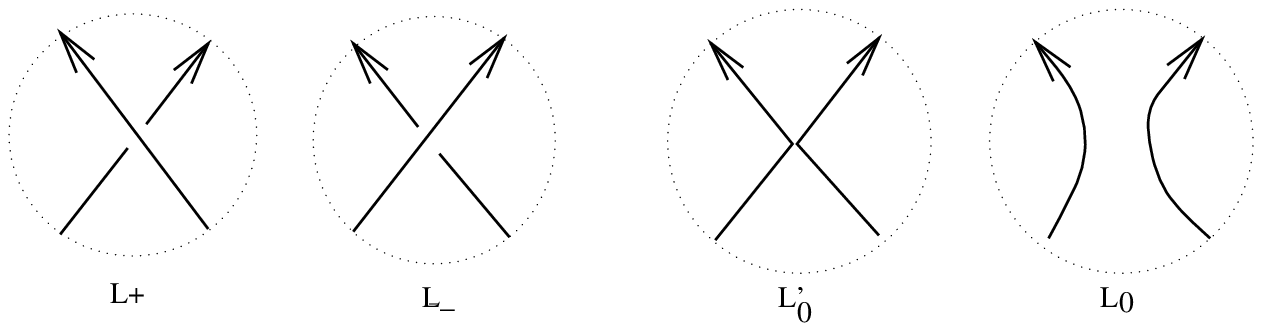}}
\par

\centerline{Fig 3}

\medskip

In addition assume that $A={\cal U}_q(sl_2)$ and that the links $L_+, L_-, L_0$
have all their components colored with the fundamental representation of $A.$
We have the
skein relations:
\be
{q}^{1\over 2}{\widehat W}_{L_{+}}-q^{-{1\over 2}}{\widehat
W}_{L_{-}}=(q-q^{-1}){\widehat W}_{L_{0}}.
\ee
Moreover:
\be
w_{L}=q^{{3\over 2}\sum_i n_i}P(L, q^2)
\ee
 where $P(L,z)$ is the Jones polynomial of the link $L$ in the variable $z.$
\ep
\proof
Let us assume first that $y_p$ and $y_q$ belongs to the same circle of the
chord diagram. We can write:
\be
\cU{1}=U_{[z_{n+1}y_{n}]}R^{(y_n^{-}y_n^{+})-1}U_{[y_n z_n]}A
U_{[z_{p+1}y_{p}]}R^{(y_p^{-}y_p^{+})-1}U_{[y_p z_p]}B.
\ee
Using this notation we obtain:
\beastar
W_{L_{+}}&=&\qTr (\prod_{i}\cU{i})\\
&=&\qTr(U_{[z_{n+1}y_{n}]}R^{(y_n^{-}y_n^{+})-1}U_{[y_n z_n]}A
U_{[z_{p+1}y_{p}]}R^{(y_p^{-}y_p^{+})-1}U_{[y_p z_p]}BC)\\
&=&\qTr(U_{[y_n z_n]}U_{[z_{n+1}y_{n}]}A
U_{[y_p z_p]}U_{[z_{p+1}y_{p}]}BC)\\
&=&\qTr(U_{[y_n z_n]}A
U_{[z_{n+1}y_{n}]}
\beta_{y_p^+}^{(y_n^{-}y_p^+)}R^{(y_n^-y_p^+)-1}
S^2(\alpha_{y_n^-}^{(y_n^{-}y_p^+)})
U_{[y_p z_p]}U_{[z_{p+1}y_{p}]}BC)\\
&=&\qTr(\beta_{y_p^+}^{(y_n^{-}y_p^+)}
U_{[y_n z_n]}A
U_{[y_p z_p]}U_{[z_{n+1}y_{n}]}
U_{[z_{p+1}y_{p}]} S^2(\alpha_{y_n^-}^{(y_n^{-}y_p^+)})BC)\\
&=&\qTr(
U_{[y_n z_n]}A
U_{[y_p z_p]}U_{[z_{n+1}y_{n}]}
U_{[z_{p+1}y_{p}]} R_{y_n^-y_p^-}^{(y_n^-y_p^+)}BC)\\
&=&\qTr(
U_{[y_n z_n]}A
U_{[z_{p+1}y_{p}]}R^{(y_p^-y_p^+)-1}U_{[y_p z_p]}
U_{[z_{n+1}y_{n}]}
R^{(y_p^-y_n^-)} R_{y_n^-y_p^-}^{(y_n^-y_p^+)}BC)\\
&=&\qTr(
U_{[y_n z_n]}A
U_{[z_{p+1}y_{p}]}R^{(y_p^-y_p^+)-1}U_{[y_p z_p]}
B R^{(z_{n+1}^-z_{n+1}^+)-1}U_{[z_{n+1}y_{n}]}
R^{(y_p^-y_n^-)} R_{y_n^-y_p^-}^{(y_n^-y_p^+)}C).
\eeastar

Using the same method, we can check that:
\beastar
W_{L_{-}}&=&\qTr(U_{[z_{n+1}y_{n}]}
R^{(y_n^+y_p^+)}R^{(y_n^{-}y_n^{+})-1}R^{(y_n^{-}y_p^{+})-1}
U_{[y_n z_n]}A\times\\
&\times& U_{[z_{p+1}y_{p}]}
R^{(y_p^{-}y_p^{+})-1}R^{(y_n^{-}y_p^{-})}R_{y_p^-y_n^-}^{(y_p^{-}y_n^{+})}
U_{[y_p z_p]}BC)\\
&=&\qTr(
U_{[y_n z_n]}A
U_{[z_{p+1}y_{p}]}R^{(y_p^-y_p^+)-1}U_{[y_p z_p]}
B R^{(z_{n+1}^-z_{n+1}^+)-1}U_{[z_{n+1}y_{n}]}
R^{(y_p^-y_n^+)}_{y_p^-y_n^-} R_{y_n^-y_p^-}^{(y_n^+y_p^+)}C).
\eeastar

{}From the standard expression of $R$ in the fundamental representation
\be
R=q^{-{1\over 2}}(q \sum_{i=1}^2 E_{ii}\otimes E_{ii}+
\sum_{i\not= j}^2 E_{ij}\otimes E_{ji} +(q-q^{-1})E_{12}\otimes E_{21})
\ee
we obtain:

\beastar
q^{1\over 2}R^{(y_p^-y_n^-)} R_{y_n^-y_p^-}^{(y_n^-y_p^+)}-
q^{-{1\over 2}}R^{(y_p^-y_n^+)}_{y_p^-y_n^-} R_{y_n^-y_p^-}^{(y_n^+y_p^+)}=
(q-q^{-1})R_{y_p^-y_n^-}^{(y_p^-y_p^+)}P_{y_n^-y_p^-}.
\eeastar
As a result we get:

\beastar
&&q^{1\over 2}W_{L_{+}}-q^{-{1\over 2}}W_{L_{-}}=\\
&=&(q-q^{-1})\qTr(P_{y_n^+y_p^+}U_{[y_n z_n]}A
U_{[z_{p+1}y_{p}]}R^{(y_p^-y_p^+)-1}U_{[y_p z_p]}BR^{(z_{n+1}^-z_{n+1}^+)}
U_{[z_{n+1}y_{n}]}R_{y_p^-y_n^-}^{(y_p^-y_p^+)}C)\\
&=&(q-q^{-1})\qTr(P_{y_n^+y_p^+}U_{[y_n z_n]}A
U_{[z_{p+1}y_{p}]}U_{[y_p z_p]}BR^{(z_{n+1}^-z_{n+1}^+)}U_{[z_{n+1}y_{n}]}C)
\eeastar
The effect of the permutation operator after taking the trace is to identify
the points $y_n^-$ with $y_p^+$ and $y_p^-$ with $y_n^+,$ therefore,

\beastar
q^{1\over 2}W_{L_{+}}-q^{-{1\over 2}}W_{L_{-}}=(q-q^{-1})W_{{L_0}'}
\eeastar

where $W_{{L_0'}}$ is the obvious generalization of the construction of Wilson
loop  to the link $L_{0}'$ (ambiant isotopic to $L_0$) whose  projection has a
non tranverse intersection at $\phi(y_p)$ (see fig 3).

We can now conclude that:
\be
{q}^{1\over 2}{\widehat W}_{L_{+}}-q^{-{1\over 2}}{\widehat
W}_{L_{-}}=(q-q^{-1}){\widehat W}_{L_{0}}.
\ee

Similar arguments would also lead to the same result when $y_p$ and $y_q$ do
not belong to the same circle.

One can easily compute the factor $v_{f}$ in the fundamental representation, we
get $v_{f}=q^{-{3\over 2}},$
as a result we obtain that $i(L)=q^{{3\over 2}\sum_i n_i}w_{L}$
is a link invariant, satisfying the skein relations $q^2 i(L_+)-q^{-2}
i(L_-)=(q-q^{-1})i(L_0).$  From uniqueness of link invariants satisfying these
axioms we get $i(L)=P(L, q^2).$
\cqfd.

\smallskip

We can now sketch a proof of the equality $w_L=RT(L)$ in the case where
$A={\cal U}_{q}(sl 2)$ for any colour of the components of $L.$
We will assume for simplicity that $L$ has only one component.
Let $\L{\alpha}\# \L{\beta}$ be the cabling of   the  framed knot $L$ with
two components coloured by $\alpha$ and $\beta.$
It can be shown (it is not completely obvious) that we have:
\be
{\widehat W}(\L{\alpha}\# \L{\beta})=
\sum_{\gamma}N_{\alpha\beta}^{\gamma}{\widehat W}(\L{\gamma})
\ee
this last equation generalizes the fusion equation (eq. \ref{fusion})  to the
case of links.
{}From the structure of the fusion ring of $sl (2)$ (it is a polynomial algebra
in one variable),  we immediately obtain that we can write:
\be
{\widehat W}(\L{\alpha})=\sum_{n=0}^{m}A_{\alpha}(n){\widehat W}(\L{f}{}^{\#
n})\label{cabling}
\ee
with $A_{\alpha}(n)\in {\bf C}.$
But we also have the same fusion equation for Reshetikhin Turaev invariant
(this is trivial from the quasitriangularity property),
as a result we also get:
\be
RT(\L{\alpha})=\sum_{n=0}^{m}A_{\alpha}(n)RT(\L{f}{}^{\# n})
\ee
with the same $A_{\alpha}(n)$ as in (eq. \ref{cabling}).
It remains to use the equality $RT(L')=w_{L'}=P(L',q^2)$ for any link with
arbitrary components coloured with the fundamental representation to obtain the
equality $RT(L)=w_L$ for any knot $L$ arbitrarily colored.
When $L$ has more than one components, it is easy to see that a generalization
of this proof works as well.

\smallskip

In the rest of this section we will give the proof of the theorem announced
at the beginning i.e in the case of the sphere the invariant of link $w_L$ is
equal to Reshetikhin-Turaev invariant. This result is quite technical and the
details are not particularly interesting in themselves.

Let $L$ be a link in $D\times [0, 1]$ where $D=[-1,1]^2$  and let choose a
braid  with $n$ colored strands which closure gives $L.$ We will denote by
$\Q{i}_{i=1,\cdots, n}$ the projections of the strands on $D$ and assume that
they are in generic positions. Let $\beta_i$ be the colour of the strand
$\Q{i}, $ these $\beta_i$ takes their values in the set $\{\alpha_j\}.$ Let
$\Delta= [0,1]\times[0,{1\over 2}] $ be the domain depicted in the picture
where all the crossings of the $\Q{i}_{i=1,\cdots, n}$ are located, and divide
$\Delta$ in $r$ strips
$\Delta_{j}=D\cap ([0,1]\times [{{r-j}\over {2r}}, {{r-j+1}\over {2r}} ]),
1\leq j\leq r$ such that in each of these strips there is only one crossing.
Let us denote by  $C_j$ and  $C_j'$  the upper horizontal part  and lower
horizontal part of the   boundary of $\Delta_j$ and by $L_j$ and  $R_j$  the
left vertical part  and right vertical part of the boundary of $\Delta_j.$
 Let us also define
$\Delta_{0}=(D\setminus\Delta)\cap\{(x,y), y\geq 0\}$ and
$\Delta_{r+1}=(D\setminus\Delta)\cap\{(x,y), y\leq 0\}.$
It will be convenient to  use the paths $\Q{i}_j=\Q{i}\cap \Delta_j,$ for $j\in
\{0,\cdots,r+1\}.$
In each of the domain $\Delta_j$ the piece of the link is an elementary braid
composed of $n$ strands with at most one crossing.

We now construct a cell decomposition of $D$ as follows:
\begin{itemize}
\item the set ${\cal V}$ of vertices is defined to be
\be
 {\cal V}=\{ x \in \partial(\Delta_i \cap \Delta_j) \,\, 0 \leq i \le j \leq
p+1 \}\cup (\cup_i \Q{i}) \cap (\cup_{j=0}^{p+1}\partial \Delta_j))\cup_{i,i'}
(\Q{i}\cap \Q{i'})
\ee

\item the set ${\cal F}$ of faces of the cell decomposition  is defined to be
the set of connected components of $D \setminus (\cup_i \Q{i}
\cup_{j=0}^{p+1}\partial \Delta_j).$
\end{itemize}

The ciliation at each vertex is shown on the following figure:

\medskip

\par
\centerline{\psfig{figure=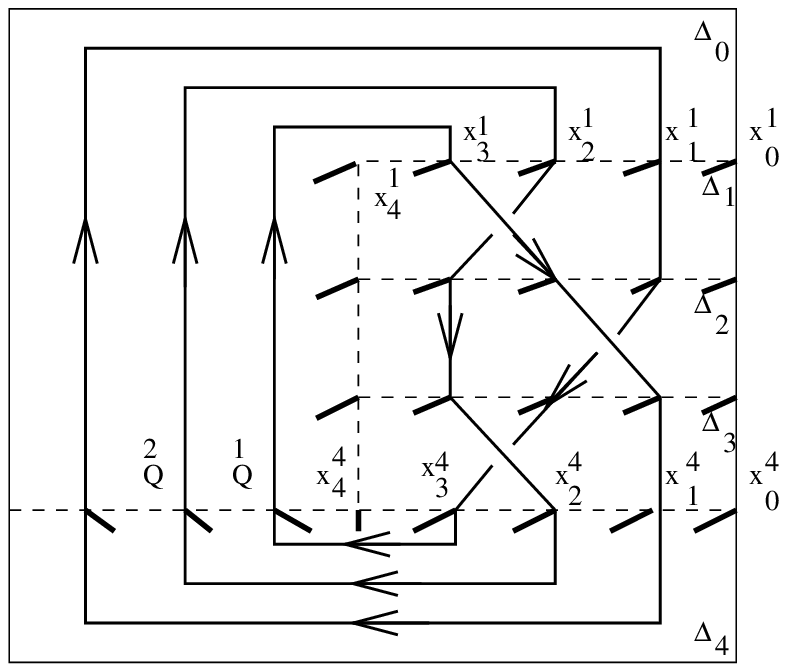}}
\par
\centerline{Fig 4}

\medskip

Let us denote by  $\beta_{i_n},\cdots, \beta_{i_1}$, in this order (from left
to right),  the incoming colors of the curves $\{\Q{i}_j, i=1\dots n\}$
according to fig 5 below, and denote by the couple $(m_j,m_j-1)$ the location
of the elementary braid (in order to simplify the notations we have put $m=m_j$
in some of the formulas below).

\medskip

\par
\centerline{\psfig{figure=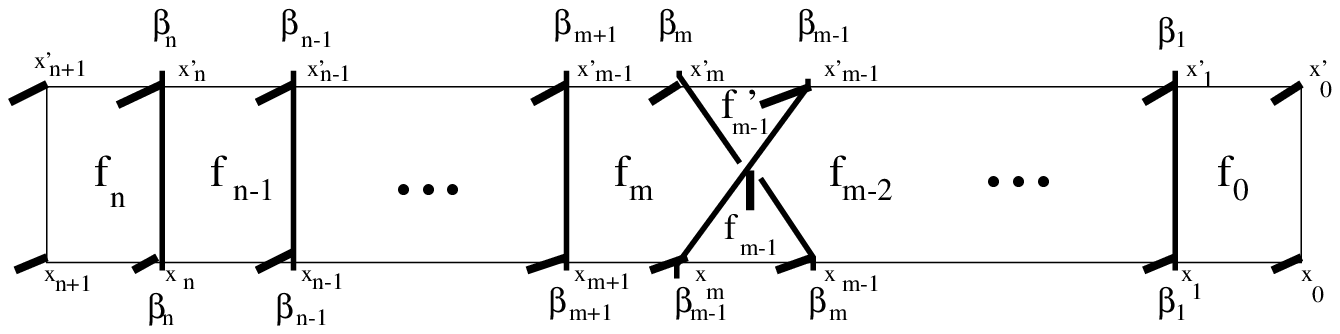}}
\par
\centerline{Fig 5}

\medskip

We can  define an element  $B_{j}$ for $0\leq j \leq r+1$
\be
B_{j}\in \Lambda\otimes End(\V{\beta_{i_1}}\otimes\cdots \otimes
\V{\beta_{i_n}},\V{\beta_{i_1}}\otimes\cdots \V{\beta_{i_{m}}}\otimes
\V{\beta_{i_{m-1}}} \otimes \cdots\V{\beta_{i_n}})
\ee

 defined by:
\bea
&&B_j=v(B_j) P_{m m-1}u(\Q{i_n}_j)_n\cdots u(\Q{i_m}_j)_m
u(\Q{i_{m-1}}_j)_{m-1}\cdots u(\Q{i_1}_j)_1\nonumber\\
&&\mbox{if $\Q{i_m}_j$ is above $\Q{i_{m-1}}_j$ and $1\leq j\leq r$,}\\
&&B_j=v(B_j) P_{m m-1}u(\Q{i_n}_j)_n\cdots u(\Q{i_{m-1}}_j)_{m-1}
u(\Q{i_m}_j)_m \cdots u(\Q{i_1}_j)_1\nonumber\\
&&\mbox{if $\Q{i_m}_j$ is under $\Q{i_{m-1}}_j$ and $1\leq j\leq r$,}\\
&&B_j=v(B_j) u(\Q{i_n}_j)_n\cdots u(\Q{i_m}_j)_m u(\Q{i_{m-1}}_j)_{m-1}\cdots
u(\Q{i_1}_j)_1\nonumber\\
&& \mbox{if $j=0$ or $r+1$}
\eea
with
\bea
&&v(B_j)=\prod_{k=1}^{n}v_{\beta_{i_k}}^{{1 \over 2}
\epsilon(d(\Q{i_k}_j),\Q{i_k})}\mbox{  if $j \not= 0$ and
}v(B_j)=\prod_{k=1}^{n}v_{\beta_{i_k}}^{-{1 \over 2}
\epsilon(d(\Q{i_k}_0),\Q{i_k})}\mbox{  if $j=0.$}\nonumber
\eea

It is not difficult to see that:

\begin{lemma}
$W_{L}$ can be expressed as:
\be
W_{L}=tr_{\V{(\beta)}}( \muf{(\beta)} \prod_{j=r+1}^0 B_{j}).
\ee
where $(\beta)=(\beta_n,\cdots, \beta_1).$
${\widehat W}_L$ can also be written:
\be
{\widehat W}_L=(\prod_{f \in {\cal F}}\delta_{\partial f}) W_{L}=
tr_{\V{(\beta)}} ( \muu{(\beta)} \prod_{j=r+1}^0 {\widehat B}_{j}),
\ee
where we have defined
${\widehat B}_{j}=
(\prod_{f \in {\cal F}\cap \Delta_j}\delta_{\partial f})B_j.$

\end{lemma}

\proof
The reader is invited to prove it. It is not difficult and uses only
the exchange relations between the edge variables $u(\Q{i}_j).$
\cqfd

\smallskip

Our strategy will consist in showing that ${\widehat W}_{L}=
(\prod_{f \in {\cal F}}\delta_{\partial f})tr_{ \V{(\beta)}}( \muu{(\beta)}
\prod_{j=r+1}^0 B_{j})$ can also be written as:
\be
{\widehat W}_{L}= (\prod_{f \in {\cal F}}\delta_{\partial f})
tr_{\V{(\beta)}}(\muu{(\beta)} \prod_{j=r+1}^0 {\Rtff{}{}}_{j})
\ee
where ${ \Rtff{}{}}_j$ is the matrix $P_{mm-1}(\beta_{m}\otimes
\beta_{m-1})(R^{\pm})$ associated to the colored braid generator defined by the
strip $\Delta_j.$
This will be sufficient to show that $w_{L}=RT(L).$
However this is not at all a trivial result and we will use integration over
the
links of the cell decomposition to obtain it.

Using again the commutation properties of proposition 1 and the specific choice
of ciliation,
we have:
\bea
&&{\widehat B}_j=\prod_{f \in {\cal F}\cap \Delta_j}\delta_{\partial
f}B_{j}=\nonumber\\
&=&v(B_j)P_{m,m-1}
(\prod_{k=n}^{m+1}\delta_{\partial f_k}u(\Q{i_k}_j)_k\delta_{\partial
f_m}\delta_{\partial {f'}_{m-1}}\delta_{\partial {f}_{m-1}}
u(\Q{i_{m-1}}_j)_{m-1} u(\Q{i_{m}}_j)_{m}\prod_{k=m-2}^{1}\delta_{\partial
f_k}u(\Q{i_k}_j)_k\nonumber
\eea

We will first establish a formula giving the expression of
\be
\int {\widehat B}_{j}
 \prod_{l\in {\cal L}\cap \{\Q{i}_j, i=1,\dots, n\}}dh (u(l))
\ee
in term of specific  elements  of the algebra $\Lambda$ (which will be
described in the sequel) and in term of the $R$ matrix expressed in the shadow
world.

Let us denote by $R^{(\pm)}_{pq}\pmatrix{\nu_3 & b' & \nu_4\cr c& & a\cr \nu_2
&b& \nu_1\cr}$, where $p, q, a, b, c, b'$ are irreducible representations of
$A$ and $\nu_1, \nu_2, \nu_3, \nu_4$ are integers labelling multiplicities, the
value of  ${\Rtff{p}{q}}{}^{\pm 1}$ in the shadow world, i.e

\be
\psi^{b'p}_{a,\nu_4}(\psi^{cq}_{b',\nu_3}\otimes id_{\V{b'}})(id_{\V{c}}\otimes
{\Rtff{p}{q}}{}^{\pm 1})(\phi^{b,\nu_2}_{cp}\otimes
id_{\V{q}})\phi^{a,\nu_1}_{bq} = R^{(\pm)}_{pq}\pmatrix{\nu_3 & b' & \nu_4\cr
c& & a\cr \nu_2 &b& \nu_1\cr} id_{\V{a}},
\ee

\begin{lemma}
The result of integration over gauge fields associated to interior links of  a
plaquette $P=A\cup B\cup C \cup B'$  describing an overcrossing:

\medskip
 \centerline{\psfig{figure=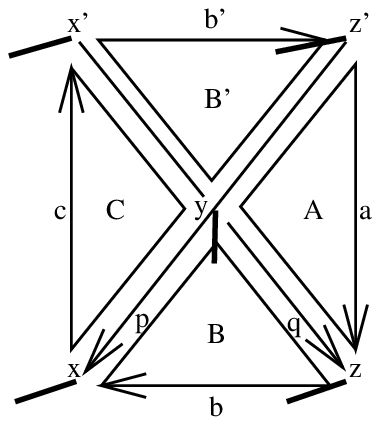}}

\centerline{Fig 6}

\medskip

has the following expression:

\begin{eqnarray}
&&\int dh(u(xy))dh(u(yz'))dh(u(zy))dh(u(yx'))\delta_{\partial A}
\delta_{\partial B}\delta_{\partial C}\delta_{\partial B'}
\u{p}(xyz')_2\u{q}(zyx')_1 \nonumber\\
&&
= {(v_p v_q)}^{-{1\over 2}}\sum_{{a,b,c,b'\in
\Irr(A)}\atop{\nu_1,\nu_2,\nu_3,\nu_4}}
[d_c]{({v_a v_{b'}\over v_{b}v_{c}})}^{1\over 2} R^{(+)}_{pq}\pmatrix{\nu_3 &
b' & \nu_4\cr c& & a\cr \nu_2 &b& \nu_1\cr}
{\cal T}_{P}\pmatrix{ \nu_3 & b' & \nu_4\cr
c& & a\cr  \nu_2&b& \nu_1\cr}\nonumber
\end{eqnarray}

where
${\cal T}_P\pmatrix{ \nu_3 & b' & \nu_4\cr
c& & a\cr  \nu_2&b& \nu_1\cr}$ is an element of $\Lambda\otimes \End
(\V{q}\otimes \V{p}, \V{q}\otimes \V{p})$ defined by:

\bea
{\cal T}_P\pmatrix{ \nu_3 & b' & \nu_4\cr
c& & a\cr  \nu_2&b& \nu_1\cr} = tr_{\V{c}}(\muf{c} \u{c}(x'x) \phi_{c
p}^{b,\nu_2} \u{b}(xz) \phi_{b q}^{a,\nu_1} \u{a}(zz')\psi^{b'p}_{a,\nu_4}
\u{b'}(z'x') \psi^{cq}_{b',\nu_1}\Rpff{c}{q}).
\eea
In the case of an undercrossing the same formula is valid after the exchange of
$\u{p}(xyz')_2\u{q}(zyx')_1$ in $\u{q}(zyx')_1\u{p}(xyz')_2$ and $R^{(+)}$ in
$R^{(-)}.$

\end{lemma}

\proof

We first show the identity:
\be
\delta_{\partial A}
\delta_{\partial B}\delta_{\partial C}\delta_{\partial B'}
(\u{p}(xyz')_1\u{q}(zyx')_2-
\u{p}(xzz')_1\u{q}(zz'x')_2)=0
\ee
which is a simple consequence of the flatness condition and the particular
choice of ciliation.

Indeed,
\bea
&&\delta_{\partial A}
\delta_{\partial B}\delta_{\partial C}\delta_{\partial B'}
\u{p}(xyz')_1\u{q}(zyx')_2=\nonumber\\
&&=\delta_{\partial B}\delta_{\partial C}\delta_{\partial B'}
\u{p}(xy)_1  v_p^{{1\over 2} \epsilon((xy)(yz))} \delta_{\partial A}
\u{p}(yzz')_1\u{q}(zyx')_2\nonumber\\
&&=\delta_{\partial A}
\delta_{\partial B}\delta_{\partial C}\delta_{\partial B'}
\u{p}(xyzz')_1\u{q}(zyx')_2\nonumber\\
&&=\delta_{\partial A}
\delta_{\partial B}\delta_{\partial C}\delta_{\partial B'}
\u{p}(xzz')_1\u{q}(zyx')_2\nonumber\\
&&=\delta_{\partial B}\delta_{\partial C}\u{p}(xzz')_1\delta_{\partial
A}\delta_{\partial B'}
\u{q}(zyx')_2\nonumber\\
&&=\delta_{\partial A}\delta_{\partial B'}
\delta_{\partial B}\delta_{\partial C}\u{p}(xzz')_1\u{q}(zz'x')_2\nonumber
\eea

Therefore we obtain:
\bea
&&\int dh(u(xy))dh(u(yz'))dh(u(zy))dh(u(yx'))\delta_{\partial A}
\delta_{\partial B}\delta_{\partial C}\delta_{\partial B'}\times\nonumber\\
&&\u{p}(xyz')_1\u{q}(zyx')_2=
\delta_{\partial P}\u{p}(xzz')_1\u{q}(zz'x')_2\nonumber
\eea
where $P$ denote the  plaquette $A\cup B\cup B'\cup C.$

\medskip

After the use of the decomposition rules of  gauge fields elements we obtain:
\begin{eqnarray*}
&&\delta_{[{x'} x z {z'}]}\guf{p}_{[x z {z'}]}\guf{q}_{[z {z'} {x'}]}{(v_p
v_q)}^{{1\over 2}}=\\
&&=\sum_{c}[d_{c}] tr_{V_{c}}(\muf{c} \guf{c}_{[{x'} {x}]}\guf{c}_{[{x}
{z}]}\guf{c}_{[{z} {z'}]}\guf{c}_{[{z'} {x'}]})
(\guf{p}_{[x z]} \guf{p}_{[z {z'}]})(\guf{q}_{[z {z'}]}v_{q}\guf{q}_{[{z'}
{x'}]})=\\
&&=\sum_{c}[d_{c}]v_{q} tr_{V_{c}}(\muf{c} \guf{c}_{[{x'} {x}]}
\guf{c}_{[{x} {z}]}\guf{p}_{[x z]}\Rpff{p}{c}\guf{c}_{[{z} {z'}]}
\guf{p}_{[z {z'}]}\guf{q}_{[z {z'}]}
\Rmff{p}{c}\Rmff{q}{c}
\guf{c}_{[{z'}{x'}]}\guf{q}_{[{z'}{x'}]})=\\
&&=\sum_{c b a b'}[d_{c}]v_{q} tr_{V_{c}}(\muf{c} \guf{c}_{[{x'} {x}]}\phi_{c
p}^{b}\guf{b}_{[{x} {z}]}
\psi_{b}^{p c}\Rtff{p}{c}\phi_{c p}^{b}\phi_{b q}^{a}\guf{a}_{[{z} {z'}]}
\psi^{q b}_{a}\Rtmff{b}{q}
\psi^{p c}_{b}
\Rtff{c}{p}\!{}^{-1} \Rpff{p}{q}
\phi_{c q}^{b'}\guf{b'}_{[{z'} {x'}]}\psi^{qc}_{b'}P_{cq})=\\
&&=\sum_{c,b,a,b'}[d_{c}] tr_{V_{c}}(\muf{c} \guf{c}_{[{x'} {x}]}\phi_{c
p}^{b}\guf{b}_{[{x} {z}]}
\phi_{b q}^{a}\guf{a}_{[{z} {z'}]}
\psi_{a}^{b' p}\guf{b'}_{[{z'} {x'}]}\psi^{cq }_{b'}\Rpff{c}{q})
R^{(+)}_{pq}\pmatrix{ & b' &\cr c& & a\cr &b& \cr}{({v_a v_{b'}\over
v_{b}v_{c}})}^{1\over 2}
\end{eqnarray*}

which establishes the result previously announced, (in order to improve the
check of this proof we have deliberately omited multiplicities).
\cqfd

\medskip

Let $\Delta_{j_2, j_1}=\cup_{j=j_1}^{j_2} \Delta_j$ for
$1\leq j_1\leq j_2\leq p,$ we can write $\partial \Delta_{j_2,j_1}=C_{j_1}\cup
R_{j_2,j_1}\cup C'_{j_2}\cup L_{j_2,j_1}$ where as usual $R_{j_2,j_1}$ and
$L_{j_2, j_1}$ denote the right vertical  part and left vertical part of the
boundary of $\Delta_{j_2,j_1}.$
It will be convenient to choose for each
$(\beta)\in Irr(A)^{\times n}, \gamma_0, \gamma_n\in Irr(A)$ a particular
basis, denoted
 $\{\psi_{\rho}, \rho\in I^{(\beta)}(\gamma_n, \gamma_0)\}$  of the spaces
$Hom_{A}(\V{\gamma_n}\otimes \V{(\beta)}, \V{\gamma_0}),$ recursively defined
by:
 $\{\psi_{\sigma}(\psi_{\gamma_n,m}^{\gamma_{n+1}\beta_{n+1}}\otimes
id_{\V{(\beta)}}),\sigma\in I^{(\beta)}(\gamma_n, \gamma_0), \gamma_{n}\in
Irr(A), m=1,.., N^{\gamma_n}_{\gamma_{n+1}\beta_{n+1}}\}$ is the basis
$\{\psi_{\rho}, \rho\in I^{(\beta_{n+1},\cdots, \beta_1)}(\gamma_{n+1},
\gamma_0)\}.$

We will also define a basis $\{\phi_{\rho}, \rho\in O^{(\beta)}(\gamma_0,
\gamma_n)\}$), of the spaces $\Hom_{A}(\V{\gamma_0}, \V{\gamma_n}\otimes
\V{(\beta)}))$ recursively defined by:
 $\{(\phi^{\gamma_n,m}_{\gamma_{n+1}\beta_{n+1}}\otimes
id_{\V{(\beta)}})\phi_{\sigma},\sigma\in O^{(\beta)}(\gamma_n, \gamma_0),
\gamma_{n}\in Irr(A), m=1,.., N^{\gamma_n}_{\gamma_{n+1}\beta_{n+1}}\}$ is the
basis
$\{\phi_{\rho}, \rho\in O^{(\beta_{n+1},\cdots, \beta_1)}(\gamma_{n+1},
\gamma_0)\}.$

If $\psi$ is an element of $Hom_{A}(\V{\gamma_n}\otimes V^{(\beta)},
\V{\gamma_0})$ and $\phi$ is an element of $Hom_{A}(\V{\gamma_n}\otimes
\V{(\beta)}, \V{\gamma_0})$ we have \be
\psi\circ \phi=a_{\psi,\phi} id_{\V{\gamma_0}}
\ee
where $a_{\psi,\phi}\in {\bf C}.$
We will use the notation $<\psi,\phi>$ to denote the number $a_{\psi,\phi}.$

The previous  choices of basis assure  that $<\psi_{\rho},\phi_{\sigma}>\in
\{0, 1\}$ if $\rho\in I^{(\beta)}(\gamma_n,\gamma_0)$ and $\sigma\in
O^{(\beta)}(\gamma_0, \gamma_n).$

We will now associate to each domain $\Delta_{j_2,j_1}$, an element of
 $\Lambda\otimes \Hom(\V{\gamma_n}\otimes \V{(\beta')}, V_{\gamma_0})$ defined
by:

\be
{\cal I}_{\Delta_{j_2,j_1}}(\rho)=
\psi_{\rho}\guf{\gamma_n}(x_{n+1}'\cdots x_{0}'x_{0})_{0}
\prod_{i=n}^{1}\guf{\beta_i'}(x_{i}\cdots x_{0}'x_{0})_{n-i+1}
\ee

and an element of $\Lambda\otimes \Hom(V_{\gamma_0},\V{\gamma_n}\otimes
\V{(\beta)})$
 defined by:

\be
{\cal O}_{\Delta_{j_2,j_1}}(\sigma)=
\prod_{i=1}^{n}\guf{\beta_i}(x_{i}\cdots
x_{0})_{i}\guf{\gamma_n}(x_{n+1}'x_{n+1}\cdots x_{0})_{0}
\phi_{\sigma}.
\ee

With these elements we can build a generalisation of the element ${\cal T}$
introduced in lemma (5) by:

\be
{\cal T}_{\Delta_{j_2, j_1}}\pmatrix{ & \rho' &\cr \gamma_{n}& & \gamma_{0}\cr
&\rho & \cr}
=tr_{\V{\gamma_n}}(\muu{\gamma_{n}} {\cal O}_{\Delta_{j_2,j_1}}(\rho){\cal
I}_{\Delta_{j_2,j_1}}(\rho'))
\ee
where $(\beta), (\beta')$ are elements of  $Irr(A)^{\times n},$ $\gamma_0,
\gamma_n$ are elements of $Irr(A)$ and $\rho\in O^{(\beta')}(\gamma_0,
\gamma_n)$ and $\rho'\in I^{(\beta)}(\gamma_n,\gamma_0)$ according to the
picture below:

\medskip

\par
\centerline{\psfig{figure=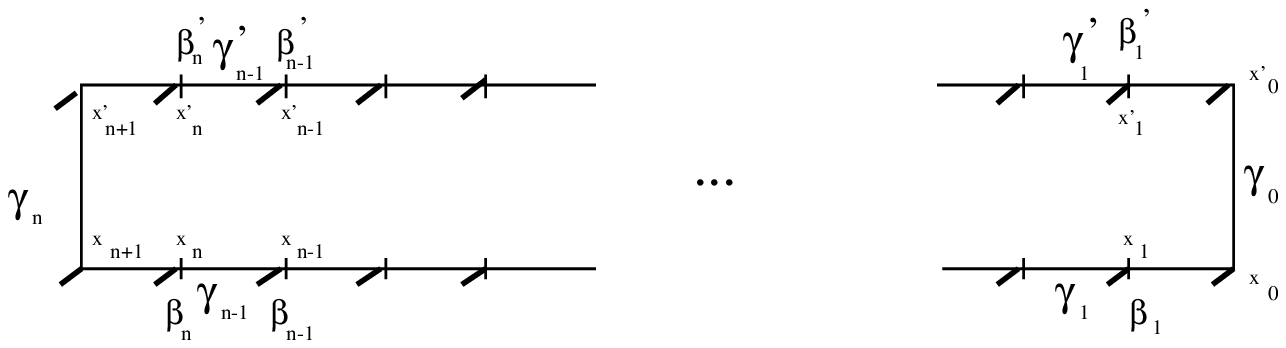}}
\par

\centerline{Fig 7}

\medskip

These elements are satisfying remarkable properties which are collected in the
next three lemmas:
\begin{lemma}
Let $j_1, j_2, j_3$ be integers such that $1\leq j_1\leq j_2\leq j_3 \leq r$
the following relation holds:

\bea
\int
&&{\cal T}_{\Delta_{j_3, j_2}}\pmatrix{ & \sigma'' &\cr \gamma_{n}'& &
\gamma_{0}'\cr &\sigma' & \cr}
{\cal T}_{\Delta_{j_2, j_1}}\pmatrix{ & \rho' &\cr \gamma_{n}& & \gamma_{0}\cr
&\rho & \cr}
\prod_{l\in {\cal L}\cap C_{j_2}} dh(u(l))=\nonumber\\
&&={1\over {[d_{\gamma_n}]}}\delta_{\gamma_0, \gamma_0'}
\delta_{\gamma_n, \gamma_n'} {\cal T}_{\Delta_{j_3, j_1}}\pmatrix{ & \sigma''
&\cr \gamma_{n}& & \gamma_{0}\cr &\rho & \cr}<\psi_{\rho'}, \phi_{\sigma'}>
\eea

where $(\beta), (\beta'), (\beta'')$ are elements of  $Irr(A)^{\times n},$
$\gamma_0, \gamma_n, \gamma_0', \gamma_n'$ are elements of $Irr(A)$ and
$\rho\in O^{(\beta)}(\gamma_0, \gamma_n), \rho'\in
I^{(\beta')}(\gamma_n,\gamma_0).$
 $\sigma'\in O^{(\beta')}(\gamma_0', \gamma_n')$ and $\sigma''\in
I^{(\beta'')}(\gamma_n',\gamma_0').$
\end{lemma}

\begin{lemma}
Consider the strip $\Delta_j$ and let us denote
 $(\beta)=(\beta_{i_n},\cdots, \beta_{i_1})$ the incoming representations.

The formula of lemma (5)  can be extended to the case of a strip $\Delta_{j}$
\be
\int \prod_{k=1}^{n}dh(u(\Q{i_k}_j))
{\widehat B}_j =
\sum_{{{\gamma_0, \gamma_n\in Irr(A)}\atop
{\rho'\in I^{(\beta)}(\gamma_n,\gamma_0)}}\atop {\rho\in
O^{(\beta)}(\gamma_0,\gamma_n)}}[d_{\gamma_n}]
{\cal T}_{\Delta_{j}}\pmatrix{ & \rho' &\cr \gamma_{n}& & \gamma_{0}\cr & \rho
& \cr}
 <\psi_{\rho'},{ \Rtff{}{}}_j \phi_{\rho}>
\ee

\end{lemma}

\medskip

\begin{lemma}[Closure of the braid]

\bea
&&\int \prod_{l \in ({\cal L} \setminus \partial D )\cap (\Delta_0 \cup
\Delta_{r+1})} dh(u_l)
 tr_{\V{(\beta)}}(\muf{(\beta)}{\widehat B}_{r+1}{\cal T}_{\Delta}\pmatrix{ &
\rho' &\cr \gamma_{n}& & \gamma_{0}\cr & \rho & \cr}
{\widehat B}_{0})= \nonumber\\
&&=\Wf{\gamma_0}_{\partial D} <\psi_{\rho'},\phi_{\rho}>
\eea
\end{lemma}

The proofs of these lemmas are not very difficult and are left to the reader.
\medskip
They imply the following proposition:

\bp[Reshetikhin-Turaev invariant]
Let $L$ be a link in $D\times[0,1],$ the element ${\widehat W}_L$ can be
 written
\be
{\widehat W}_L=w_L \prod_{f\in {\cal F}}\delta_{\partial f}
\ee
where $w_{L}=RT(L).$
\ep
\proof
{}From lemmas (6,7) we obtain:
\be
\int \prod_{l\in (\Delta\setminus \partial \Delta)\cap {\cal L}}
dh(u_l)\prod_{j=r}^{1}{{\widehat B}_j}=
\sum_{{{\gamma_0, \gamma_n\in Irr(A)}\atop
{\rho'\in I^{(\beta)}(\gamma_n,\gamma_0)}}\atop {\rho\in
O^{(\beta)}(\gamma_0,\gamma_n)}}
{\cal T}_{\Delta}\pmatrix{ & \rho' &\cr \gamma_{n}& & \gamma_{0}\cr & \rho &
\cr}
 <\psi_{\rho'},\prod_{j=r}^{1}{ \Rtff{}{}}_j\phi_{\rho}>.
\ee
Using lemma (8) we get:

\be
\int \prod_{l\in (D\setminus \partial D)\cap {\cal L}}
dh(u_l)tr_{\V{(\beta)}}tr(\muf{(\beta)}\prod_{j=r+1}^{0}{{\widehat B}_j})=
\sum_{
{{\gamma_0, \gamma_n\in Irr(A)}\atop
{\rho'\in I^{(\beta)}(\gamma_n,\gamma_0)}}\atop
 {\rho\in O^{(\beta)}(\gamma_0,\gamma_n)}
} [d_{\gamma_n}]
\Wf{\gamma_0}_{\partial D} <\psi_{\rho'},\phi_{\rho}>
 <\psi_{\rho'},\prod_{j=r}^{1}{ \Rtff{}{}}_j\phi_{\rho}>.
\ee
{}From the relation:

\be
\sum_{
{{ \gamma_n\in Irr(A)}\atop
{\rho'\in I^{(\beta)}(\gamma_n,\gamma_0)}}\atop {\rho\in
O^{(\beta)}(\gamma_0,\gamma_n)}
}
{[d_{\gamma_n}]}
<\psi_{\rho'},\phi_{\rho}>\phi_{\rho}(id_{\V{\gamma_n}}\otimes A)\psi_{\rho'}
=
[d_{\gamma_0}]id_{\V{\gamma_0}}tr_{\V{(\beta)}}(\muf{(\beta)}A)
\ee

which holds for every $A\in End(\V{(\beta)})$ (this relation is a consequence
of relations (13,14) on the Clebsh-Gordan maps (see proposition 14 of
\cite{BR}))
we obtain:
\be
\int \prod_{l\in (D\setminus \partial D)\cap {\cal L}}
dh(u_l)tr_{\V{(\beta)}}tr(\muf{(\beta)}\prod_{j=r+1}^{0}{{\widehat B}_j})=
(\sum_{\gamma_0} [d_{\gamma_0}]\Wf{\gamma_0}_{\partial
D})tr_{\V{(\beta)}}(\muf{(\beta)}\prod_{j=r}^{1}{ \Rtff{}{}}_j)
\ee
which is the last step in the proof.
\cqfd

\section{Conclusion}
In this work we continued the  analysis of  combinatorial quantization of
hamiltonian Chern-Simons theory. We have defined elements of the observable
algebra $\Lambda_{CS}$ which are associated to any link in $\Sigma\times {\bf
R}.$ These elements are expected to be the precise definition of the Wilson
loop elements in E.Witten formalism. This is supported by the fact that they
are ribbon invariants and that their expectation value when $\Sigma=S^2$ is
precisely the Reshetikhin-Turaev invariant of the link.
In the next work \cite{Bu} we have  studied this combinatorial approach when
the 3-manifold $M$ is arbitrary. It will be much more convenient to work with a
Heegaard splitting of M rather than using a surgery presentation. If $\Sigma_g$
and $f\in Mod(\Sigma_g)$  is any Heegaard splitting of $M$ we can associate an
element $a_f$ of $\Lambda_{CS}(\Sigma_g)$ which expectation value  gives an
invariant of $M$, which can be shown to be the Reshetikhin -Turaev invariant of
$M$. This invariant can be expressed as the partition function of a non
commutative lattice gauge theory associated to   a cellular decomposition of
$M$. One important technical point which we will have to develop is the
truncation of the spectrum when $q$ is a root of unity.

The expression of the observables associated to any link that we have found  is
derived from few principles: gauge invariance, independance under the choice of
departure points. We expect that similar arguments will lead to the
construction of observables associated to links with arbitrary self crossings.
This could shed lights on  the construction of states  in the canonical
quantization program of pure gravity \cite{Ba}.

\medskip

{Acknowledgments}

We warmfully thank N.Reshetikhin for fruitful discussions and suggestions.
We also  thank C.Mercat for helpful remarks on the first draft of this work.

\newpage

\bibliographystyle{unsrt}

\begin{thebibliography}{10}


\bibitem{AMR1}
J.E.Andersen, J.Mattes, N.Reshetikhin,
\newblock{The Poisson Structure on the Moduli Space of Flat Connections and
Chord Diagrams }
\newblock{To appear} {\bf }, (1995).


\bibitem{AMR2}
J.E.Andersen, J.Mattes, N.Reshetikhin,
\newblock{A Quantization of the Moduli Space of Flat Bundles over Surfaces }
\newblock{Unpublished} {\bf }, (1995).







\bibitem{AGS1}
A.Y.Alekseev, H.Grosse, V.Schomerus,
\newblock{Combinatorial Quantization of the Hamiltonian Chern-Simons Theory I,}
\newblock{\it hep-th /94/03,} {\bf } (1994).

\bibitem{AGS2}
A.Y.Alekseev, H.Grosse, V.Schomerus,
\newblock{Combinatorial Quantization of the Hamiltonian Chern-Simons Theory
II,}
\newblock{\it hep-th /94/08,} {\bf } (1994).

\bibitem{AS}
A.Y.Alekseev,  V.Schomerus,
\newblock{Representation Theory of Chern-Simons Observables,}
\newblock{\it hep-th /95/03,} {\bf } (1995).

\bibitem{AxS}
S. Axelrod, I. Singer,
\newblock{Chern-Simons perturbation theory},
\newblock{Proceedings of the XXth International conference on differential
geometric methods in theoretical physics, June 3-7} New York City, (1991)
(World Scientific, Singapore, 1992).

\bibitem{Ba}
J.Baez, Editor,
\newblock{Knots and Quantum Gravity,}
\newblock{Oxford Lectures Series in Mathematics and its Applications
n${}^{\circ}$ 1.}




\bibitem{Bu}
E.Buffenoir,
\newblock{Chern-Simons Theory on a Lattice as a new description of 3-Manifold
Invariants},
\newblock{To appear}.




\bibitem{BN}
D.Bar-Natan,
\newblock{On the Vassiliev knot invariants,}
\newblock{\it Topology,} {\bf 34} $N^{\circ}2$, (1995).





\bibitem{BR}
E.Buffenoir, Ph.Roche,
\newblock{Two dimensional lattice gauge theory based on a quantum group,}
\newblock{\it Comm.Math.Phys,} {\bf 170} (1995).




\bibitem{Dr}
V.G.Drinfeld,
\newblock{On almost cocomutative Hopf algebras,}
\newblock{\it Leningrad.Math.Journal,} {\bf Vol1. nb 2}, 321 (1990).

\bibitem{ES}
D.Eliezer, G.W.Semenoff,
\newblock{Intersection forms and the geometry of lattice Chern-Simons theory,}
\newblock{Phys.Lett.B,} {\bf 286}, 118-124 (1992)



\bibitem{FR}
V.V.Fock, A.A.Rosly,
\newblock{Poisson structure on moduli of flat connections on Riemann surfaces
and
r-matrices,} \newblock{\it Preprint ITEP 72-92,} {\bf }, (1992).


\bibitem{GMM}
E. Guadagnini, M. Martellini and M. Mintchev,
\newblock{Perturbative aspect of Chern-Simons field Theory,}
\newblock{\it Phys. Lett.} {\bf B227} (1989)





\bibitem{Ma}
J.M.Maillet,
\newblock{Lax equations and quantum groups}
\newblock{\it Phys.Lett.B,}{\bf 245}, (1990)




\bibitem{RT1}
N.Yu.Reshetikhin, V.G.Turaev,
\newblock{Ribbon Graphs and their invariant derived from quantum groups,}
\newblock{\it Comm.Math.Phys,} {\bf 127}, (1990).

\bibitem{RT2}
N.Yu.Reshetikhin, V.G.Turaev,
\newblock{Invariants of 3-manifolds via link polynomials and quantum groups,}
\newblock{Invent.Math.} {\bf 103},  547-597 (1991).


\bibitem{T}
 V.G.Turaev,
\newblock{Skein quantization of Poisson algebras of loops on surfaces,}
\newblock{Ann.Scient.ENS} {\bf 24},  635-704 (1991).

\bibitem{TV}
 V.G.Turaev, O.Y.Viro
\newblock{State Sum Invariants of 3-Manifolds and quantum 6j-symbols,}
\newblock{Topology} {\bf 31}, 865-902 (1992).



\bibitem{Va}
V.A.Vassiliev,
\newblock{Topology of Complements to Discriminants and Loop Spaces,}
\newblock{\it Adv.Sov.Math,} {\bf 1}, 9 (1990).


\bibitem{W}
E.Witten,\newblock{Quantum field theory and Jones polynomial,}
\newblock{\it Comm.Math.Phys,} {\bf 121}, 351-399 (1989).














\end{thebibliography}

\end{document}